\documentclass[10pt]{article}

\date{}

\usepackage{latexsym,amssymb}
\usepackage{amsmath}
\usepackage{mathrsfs}
\usepackage{euscript}
\usepackage{bussproofs}
\usepackage{tikz}
\usepackage{amsfonts}   

\usepackage[bottom]{footmisc}

\newtheorem{defi}{Definition}[section]
\newtheorem{theo}[defi]{Theorem}
\newtheorem{lem}[defi]{Lemma}
\newtheorem{prop}[defi]{Proposition}
\newtheorem{rem}[defi]{Remark}

\newtheorem{cor}[defi]{Corollary}

\def\Sup{\mathop{\rm Sup}\nolimits}
\def\Inf{\mathop{\rm Inf}\nolimits}

\newcommand{\negr}[1]{\boldsymbol{#1}}

\newenvironment{dem}{\noindent\bf Proof. \rm}{\hfill $\negr{\blacksquare}$}

 \relax
\input{prepictex}
\input{pictex}
\input{postpictex}

\title{Cut--free sequent calculus and natural deduction for the tetravalent modal logic}

\author{Mart\'in Figallo}

\date{\textit{\small Departamento de Matem\'atica. Universidad Nacional del Sur. Bah\'{\i}a Blanca, Argentina}}

\begin{document}

\maketitle

\begin{abstract}
The {\em tetravalent modal logic} ($\cal TML$) is one of the two logics defined by Font and Rius (\cite{FR2}) (the other is the {\em normal tetravalent modal logic} ${\cal TML}^N$) in connection with Monteiro's tetravalent modal algebras. These logics are expansions of the well--known {\em Belnap--Dunn's four--valued logic} that combine a many-valued character (tetravalence)  with a modal character. In fact, $\cal TML$ is the logic that preserve degrees of truth with respect to tetravalent modal algebras. As Font and Rius observed, the connection between the logic $\cal TML$ and the algebras is not so good as in ${\cal TML}^N$, but, as a compensation, it has a better proof-theoretic behavior, since it has a strongly adequate Gentzen calculus (see \cite{FR2}). In this work, we prove that the sequent calculus given by Font and Rius does not enjoy the cut--elimination property. Then, using a general method proposed by Avron, Ben-Naim and Konikowska (\cite{Avron02}), we provide a sequent calculus for $\cal TML$ with the cut--elimination property. Finally, inspired by the latter, we present a {\em natural deduction} system, sound and complete with respect to the tetravalent modal logic.
\end{abstract}

\

\section{Introduction}\label{s1}\indent

The class {\bf TMA} of tetravalent modal algebras was first
considered by Antonio Monteiro (1978), and mainly studied by I.
Loureiro, A.V. Figallo, A. Ziliani and P. Landini. Later on, J.M.
Font and M. Rius were interested in the logics arising from the
algebraic and lattice--theoretical aspects of these algebras.
From Monteiro's point of view,  in the future these algebras 
would give rise to a four-valued modal logic with significant applications in Computer Science (see \cite{FR2}).
Although such applications have not yet been developed, the two logics considered in \cite{FR2} are modal expansions
of Belnap-Dunn's four-valued logic, a logical system that is well--known for the many applications it has found in several fields. 
In these logics, the four non-classical epistemic values emerge: 1 (true and not false), 0 (false and not true), {\bf n} (neither true nor false) and {\bf  b} (both true and false). 
We may think of them as the four possible ways in which an atomic sentence $P$ can
belong to the {\em present state of information} : we were told that (1)  $P$ is true (and were not told that $P$ is false); (2) $P$
is false (and were not told that $P$ is true); (3) $P$ is both true and false (perhaps from different sources, or in different instants
of time); (4) we were not told anything about the truth value of $P$.
In this interpretation, it makes sense to consider a modal-like unary operator $\square$ of epistemic character, such that for any sentence $P$, the sentence $\square P$ would mean  ``the available information confirms that $P$ is true".
 It is clear that in this setting the sentence $\square P$ can only be true in the case where we have some information saying that $P$  is true and we have no information saying that $P$ is false, while it is simply false in all other cases (i.e., lack of information or at least some information saying that $P$ is false, disregarding whether at the same time some other information says that $P$ is true); that is, on the set $\{0, {\bf n}, {\bf b}, 1\}$ of epistemic values this operator must be defined as $\square 1 = 1$ and $\square {\bf n} = \square {\bf b} =\square 0=0$ . This is exactly the algebra that generates the variety of TMAs. 

In \cite{FR2}, Font and Rius studied two logics related to TMAs. One of them is obtained by following the usual ``preserving truth''  scheme,  taking  $\{1\}$ as designated  set, that is, $\psi$ follows from $\psi_1, \dots, \psi_n$ in  this  logic  when every interpretation  that  sends  all  the  $\psi_i$ to  $1$  also sends  $\psi$ to  $1$. The  other logic, denoted by ${\cal TML}$ (the logic we are interested in), is defined by  using  the  {\em preserving  degrees of truth}  scheme, that is, $\psi$  follows  from  $\psi_1, \dots, \psi_n$  when  every interpretation that assigns to  $\psi$ a value that is greater  or equal than the value it assigns to  the conjunction of  the  $\psi_i$'s.  These authors proved that ${\cal TML}$  is  not  algebraizable  in  the  sense of Blok and Pigozzi, but it is  {\em finitely equivalential}  and  {\em protoalgebraic}. However, they confirm that its algebraic  counterpart  is also the  class of TMAs: but the  connection between the logic and the algebras is not so good  as in the first logic.  As a compensation,  this logic has a better  proof-theoretic  behavior,  since  it  has  a {\em strongly  adequate  Gentzen 
calculus} (Theorems  3.6 and 3.19 of \cite{FR2}). 

In \cite{FR2}, it was proved that $\cal TML$ can be characterized as a matrix logic in terms of two logical matrices, but later, in \cite{CF}, it was proved that $\cal TML$ can be determined by a single logical matrix. Besides, taking profit of the contrapositive implication introduced by A. V. Figallo and P. Landini (\cite{FL}), a  sound and complete Hilbert-style calculus for this logic was presented. Finally, the paraconsistent character of  $\cal TML$ was also studied from the point of view of the \emph{Logics
of Formal Inconsistency}, introduced by W. Carnielli and J. Marcos in~\cite{Tax} and afterward developed in~\cite{WCMCJM}. 

\section{Preliminaries}\label{s11}\indent

Recall that, a {\em De Morgan algebra} is a structure $\langle
A, \wedge, \vee, \neg, 0 \rangle$ such that $\langle A, \wedge, \vee, 0\rangle$ 
is a bounded distributive lattice and
$\neg$ is a De Morgan negation, i.e., an involution that
additionally satisfies De Morgan's laws: for every $a, b\in A$
$$\neg \neg a = a$$ 
$$\neg (a \vee b) = \neg a \wedge  \neg b.$$

A {\em tetravalent modal algebra} (TMA) is an algebra
$\mathbb{A}=\langle A, \wedge, \vee, \neg, \square, 0\rangle$ of
type $(2,2,1,1,0)$ such that its non-modal reduct
$\langle A, \wedge, \vee, \neg, 0\rangle$ is a De
Morgan algebra and the unary operation $\square$ satisfies, for
all $a\in A$, the two following axioms:
 $$\square a \wedge \neg a = 0,$$ 
 $$\neg \square a \wedge a = \neg a \wedge a.$$ 

Every TMA $\mathbb{A}$ has a top element $1$ which is defined as $\neg 0$. 
These algebras were studied mainly by I. Loureiro (\cite{L}), and also by A. V.  Figallo, P. Landini (\cite{FL}) and A. Ziliani, at the suggestion of the late A. Monteiro (see \cite{FR2}).
The class of all tetravalent modal algebras constitute
a variety which is denoted by {\bf TMA}. Let $M_4=\{0,{\bf n}, {\bf b}, 1\}$ and consider the lattice given by the following Hasse diagram 
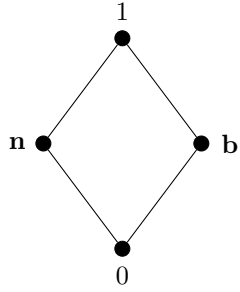
\begin{center}
\begin{tikzpicture}[scale=.7]
\tikzstyle{every node}=[draw,circle,fill=black,inner sep=2pt]
  \node (one) at (0,4) [label=above:$1$] {};
  \node (b) at (-1.5,2) [label=left:$\bf n$] {};
   \node (a) at (1.5,2) [label=right:$\bf b$] {};
  \node (zero) at (0,0) [label=below:$ 0$] {};
  \draw (zero) -- (a) -- (one) -- (b) -- (zero);
\end{tikzpicture} \hspace{2cm}
\end{center}
This is a well-known lattice and it is called ${\bf L4}$ (See \cite{AnBel}, pg. 516.) 
Then, {\bf TMA} is generated by the above  four--element   lattice enriched with two unary operators $\neg$ and $\square$
 given by $\neg {\bf n}= {\bf n}$, $\neg {\bf b}={\bf b}$, $\neg 0=1$ and $\neg 1=0$ and the unary operator $\square$ is defined as: $\square {\bf n}= \square {\bf b}= \square 0=0$ and $\square 1=1$  (see \cite{FR2}). This tetravalent modal algebra, denoted by  $\mathfrak{M}_{4m}$, has two prime filters, namely, $F_{\tiny \mbox{\bf n}} =\{{\bf n}, 1\}$ and
$F_{\tiny \mbox{\bf b}} =\{{\bf b}, 1\}$. As we said, $\mathfrak{M}_{4m}$ generates the variety ${\bf TMA}$, i.e., an
equation holds in every  TMA  iff it holds in $\mathfrak{M}_{4m}$. 

\begin{lem}(See \cite{FR2})\label{lem2.1}\label{PropTMA} In every TMA $\mathbb{A}$ and for all $a,b\in A$ the following hold:\\[2mm]
\begin{tabular}{clcl}
{\bf (i)} &  $\neg \square a \vee  a=1$, \hspace{1cm} & {\bf (viii)} & $\square \square a =  \square a$,\\
{\bf (ii)} & $\square a \vee  \neg a = a \vee  \neg a $,  & {\bf (ix)} & $\square(a\wedge b)=\square a\wedge \square b$, \\
{\bf (iii)} & $\square a \vee \neg\square a=1$, & {\bf (x)} & $\square (a \vee \square b) = \square a \vee \square b$,\\
{\bf (iv)} & $\square a \wedge \neg \square a =0$, & {\bf (xi)} & $\square \neg \square a = \neg\square a$ \\
{\bf (v)} &  $\square a \leq  a$, & {\bf (xii)} & $ a\wedge \square \neg a = 0$,\\
{\bf (vi)} & $\square 1 = 1$,  & {\bf (xiii)} & $\square(\square a\wedge \square b)=\square a\wedge \square b$ \\
{\bf (vii)} & $\square 0 = 0$, & {\bf (xiv)} & $\square(\square a\vee \square b)=\square a\vee \square b$ \\
\end{tabular}
\end{lem}

The next proposition will be needed in what follows.

\newpage

\begin{prop} \label{propsquare} Let $\mathbb{A}$ be a TMA. If $x\leq y\vee z$ and $x\wedge \neg z\leq y$, then $x\leq y \vee \square z$, for every $x,y,z\in A$.
\end{prop}
\begin{dem} It is a routine task to check that the assertion holds in $\mathfrak{M}_{4m}$. The fact that $\mathfrak{M}_{4m}$ generates the variety ${\bf TMA}$ completes the proof.
\end{dem}
 
\

Let $\mathscr{L}=\{\vee, \wedge, \neg, \square\}$ be a propositional language. 
From now on, we shall denote by $\mathfrak{Fm}=\langle Fm, \wedge,
\vee, \neg, \square, \bot \rangle$ the absolutely free algebra of
type (2,2,1,1,0) generated by some denumerable set of variables.
We denote by $Fm$ the set of sentential formulas, and we shall refer to them
by lowercase Greek letters $\alpha, \beta, \gamma, \dots$ and so
on; and we shall denote finite sets of formulas by uppercase Greek
letters $\Gamma, \Delta,$ etc.

\begin{defi} \label{logTMA} The tetravalent modal logic ${\cal TML}$
defined over $\mathfrak{Fm}$ is the propositional logic
$\langle Fm, \models_{{\cal TML}}\rangle$ given as
follows: for every finite set $\Gamma \cup\{\alpha\} \subseteq
Fm$, $\Gamma \models_{\cal TML} \alpha$  if and only if, for every
$\mathbb{A} \in {\bf TMA}$ and for every $h \in
Hom(\mathfrak{Fm}, \mathbb{A})$, $\bigwedge \{ h(\gamma) \ : \
\gamma \in \Gamma\} \leq h(\alpha)$. In particular, $\emptyset
\models_{\cal TML} \alpha$ if and only if $h(\alpha)=1$ for every
$\mathbb{A} \in {\bf TMA}$ and for every $h \in
Hom(\mathfrak{Fm}, \mathbb{A})$.
\end{defi}

\begin{rem} Observe that, if $h\in Hom(\mathfrak{Fm}, \mathbb{A})$ for any $\mathbb{A} \in {\bf TMA}$, we have that $h(\bot)=0$. This follows from the fact that $\bot$ is the $0$-ary operation in $\mathfrak{Fm}$, $0$ is the $0$-ary operation in  $\mathbb{A}$ and the definition of homomorphism (in the sense of universal algebra). 
\end{rem}

\

Let ${\cal M}=\langle {\cal T}, {\cal D }, {\cal O}\rangle$ be a logical matrix for $\mathscr{L}$, that is, ${\cal T}$ is a finite, non-empty set of truth values, ${\cal D}$ is a non-empty proper subset of ${\cal T}$, and ${\cal O}$ includes a $k$-ary function $\hat{f}: {\cal T}^k\to {\cal T}$ for each $k$-ary connective $f\in \mathscr{L}$. Recall that, a valuation in ${\cal M}$ is a function $v:Fm\to {\cal T}$ such that 
$$v(f(\psi_1, \dots, \psi_k))=\hat{f}(v(\psi_1), \dots, v(\psi_k))$$
for each $k$-ary connective $f$ and all $\psi_1, \dots, \psi_k\in Fm$. A formula $\alpha \in Fm$ is satisfied by a given valuation $v$, in symbols $v\models \alpha$, if $v(\alpha)\in {\cal D}$. Let $\Gamma, \Delta \subseteq Fm$. We  say that the  $\Delta$ is consequence of $\Gamma$, denoted $\Gamma \models_{\cal M} \Delta$, iff  for every valuation $v$ in ${\cal M}$, either $v$ does not satisfy some formula in $\Gamma$ or $v$ satisfies some formula in $\Delta$. \\[2mm]
J. M. Font and M. Rius proved in \cite{FR2} that the tetravalent modal logic $\cal TML$ is a matrix logic defined in
terms of two logical matrices. But later, M. E. Coniglio and M. Figallo  proved in  \cite{CF} that $\cal TML$ can be characterized as a matrix logic in terms of a single logical matrix. Indeed, let ${\cal M}_4=\langle {\cal T}, {\cal D}, {\cal O} \rangle$ be the matrix where the set of truth values is  ${\cal T}=\{0,{\bf n}, {\bf b}, 1\}$, the set of designated values is ${\cal D}=\{{\bf b}, 1\}$ and ${\cal O}=\{\tilde{\vee}, \tilde{\wedge}, \tilde{\neg}, \tilde{\square}\}$ where $\tilde{\vee}, \tilde{\wedge}:{\cal T}^2\to{\cal T}$ and $\tilde{\neg}, \tilde{\square}:{\cal T}\to{\cal T}$ are defined as $x\tilde{\vee}y=\Sup\{x,y\}$, $x\tilde{\wedge}y=\Inf\{x,y\}$ (here we are assuming that the elements of ${\cal T}$ are ordered as in the lattice $M_4$).   

\begin{center}
\begin{tabular}{c|c|c} 
$x$ & $\tilde{\neg} x$ & $\tilde{\square} x$  \\ \hline
0 & 1 & 0\\ 
{\bf n} & {\bf n} &0\\ 
{\bf b} & {\bf b} &0\\ 
1 & 0 &1\\ 
\end{tabular}
\end{center}

\

\noindent then, 

\begin{prop}\label{PropCompTML} (\cite{CF}) $\cal TML$ is sound and complete w.r.t. ${\cal M}_4$.
\end{prop}

Therefore, given $\Gamma$ and $\Delta$ sets of formulas,  {\em $\Delta$ is consequence of $\Gamma$ in ${\cal TML}$}, denoted $\Gamma \models_{\cal TML} \Delta$, iff  for every valuation $v$ in ${\cal M}_4$, either $v$ does not satisfy some formula in $\Gamma$ or $v$ satisfies some formula in $\Delta$. If $\Delta$ is a set with exactly one element, we recover the consequence relation given in Definition \ref{logTMA}.
 
\

In order to characterize $\cal TML$ syntactically, that is, by means of a deductive system, J. M. Font and M. Rius introduced in~\cite{FR2} the sequent calculus $\mathfrak{G}$. The sequent calculus $\mathfrak{G}$  is single--conclusion, that is, it deals with sequents of the form $\Delta \Rightarrow \alpha$ such that $\Delta \cup \{\alpha\}$ is a finite subset of $Fm$. The axioms and rules of $\mathfrak{G}$ are the following:

\

\noindent {\bf Axioms}
$$ \mbox{(Structural axiom) \, } \displaystyle {\alpha \Rightarrow \alpha} \hspace{2cm} \mbox{(Modal axiom) \,  } {\Rightarrow \alpha \vee \neg \square \alpha}$$

\noindent {\bf Structural rules}

$$ \mbox{(Weakening) \, } \displaystyle \frac{\Delta \Rightarrow \alpha} {\Delta, \beta \Rightarrow \alpha} \hspace{2cm} \mbox{(Cut) \, }
 \displaystyle \frac{\Delta \Rightarrow \alpha  \hspace{0.5cm} \Delta, \alpha \Rightarrow \beta}{\Delta \Rightarrow \beta} $$

\noindent {\bf Logic rules}

$$ \mbox{($\wedge \Rightarrow$) \, } \displaystyle \frac{\Delta, \alpha, \beta \Rightarrow \gamma} {\Delta, \alpha \wedge \beta \Rightarrow \gamma} \hspace{2cm} \mbox{($\Rightarrow \wedge$) \, }
 \displaystyle \frac{\Delta \Rightarrow \alpha  \hspace{0.5cm} \Delta \Rightarrow \beta}{\Delta \Rightarrow \alpha \wedge \beta} $$

$$ \mbox{($\vee \Rightarrow$) \, } \displaystyle \frac{\Delta, \alpha \Rightarrow \gamma \hspace{0.5cm} \Delta, \beta \Rightarrow \gamma} {\Delta, \alpha \vee \beta \Rightarrow \gamma}$$

$$ \mbox{($\Rightarrow \vee$)$_1$ \, }
 \displaystyle \frac{\Delta \Rightarrow \alpha }{\Delta \Rightarrow \alpha \vee \beta} \hspace{2cm} \mbox{($\Rightarrow \vee$)$_2$ \, } \displaystyle \frac{\Delta \Rightarrow \beta }{\Delta \Rightarrow \alpha \vee \beta}$$

$$ \mbox{($\neg$) \, } \displaystyle \frac{\alpha \Rightarrow \beta } {\neg \beta \Rightarrow \neg \alpha }
 \hspace{2cm}  \mbox{($\bot$) \,} \frac{\Delta \Rightarrow \bot } {\Delta \Rightarrow \alpha}$$

$$ \mbox{($\neg \neg \Rightarrow$) \, } \displaystyle \frac{\Delta, \alpha \Rightarrow \beta } {\Delta, \neg \neg \alpha \Rightarrow \beta} \hspace{2cm} \mbox{($\Rightarrow \neg \neg$)} \, \frac{\Delta \Rightarrow \alpha } {\Delta \Rightarrow \neg \neg \alpha} $$

$$ \mbox{($\square \Rightarrow$) \, } \displaystyle \frac{\Delta, \alpha, \neg \alpha \Rightarrow \beta } {\Delta, \alpha, \neg \square \alpha \Rightarrow \beta} \hspace{2cm} \mbox{($\Rightarrow \square$)} \, \frac{\Delta \Rightarrow \alpha \wedge \neg \alpha } {\Delta \Rightarrow \alpha \wedge \neg \square \alpha} $$

\

\noindent The notion of derivation in  the sequent calculus $\mathfrak{G}$ is the usual. Besides, for every finite set $\Gamma \cup
\{\varphi\} \subseteq Fm$, we write $\Gamma \vdash_{\mathfrak{G}}  \varphi$ iff the sequent $\Gamma \Rightarrow \varphi$ has a derivation in $\mathfrak{G}$. We say that the sequent $\Gamma \Rightarrow \varphi$ is provable iff there exists a derivation for it in $\mathfrak{G}$.\\[2mm]
J. M. Font and M. Rius proved in~\cite{FR2} that  $\mathfrak{G}$ is sound and complete with respect to the tetravalent modal logic $\cal TML$.

\

\begin{theo}\label{compTML} {\rm(Soundness and Completeness, \cite{FR2})} For every finite set $\Gamma
\cup\{\alpha\} \subseteq Fm$,
$$\Gamma \models_{\cal TML} \alpha \ \ \textrm{ if and only if} \ \ \Gamma \vdash_{\mathfrak{G}} \alpha.$$
\end{theo}

\

\noindent Moreover,

\begin{prop}{\rm(\cite{FR2})}
An arbitrary equation $\psi \approx \varphi$ holds in every TMA iff $\psi \dashv\vdash_{\mathfrak{G}}\varphi$ (that is, $\psi \vdash_{\mathfrak{G}}\varphi$ and $\varphi \vdash_{\mathfrak{G}}\psi$).
\end{prop}

\noindent As a consequence of it we have that:

\begin{cor}{\rm(\cite{FR2})}\label{CorTheoiffValid}
\begin{itemize}
\item[{\rm(i)}] The equation $\psi \approx 1$ holds in every TMA iff \, $\vdash_{\mathfrak{G}}\psi$.
\item[{\rm(ii)}] For any $\psi, \varphi \in Fm$,\,  $\psi \vdash_{\mathfrak{G}} \varphi$ \, iff \,
 $h(\psi) \leq h(\varphi)$ \, for every \, $h \in Hom(\mathfrak{Fm}, \mathbb{A})$, \, for every
 $\mathbb{A} \in {\bf TMA}$.
\end{itemize}
\end{cor}

\

\section{$\mathfrak{G}$ does not admit a cut--elimination theorem}\label{s2}

Corollary~\ref{CorTheoiffValid} is a powerful tool to determine whether a given sequent of $\mathfrak{G}$ is provable or not. For instance,

\begin{prop}\label{prop1} In $\mathfrak{G}$ we have that the sequent $\neg \square \alpha \Rightarrow \alpha$ is provable iff \, the sequent $\Rightarrow \alpha$ is provable.
\end{prop}
\begin{dem}
 Indeed, suppose that the sequent $\neg \square \alpha \Rightarrow \alpha$ is provable in $\mathfrak{G}$. Then, \, $h(\neg \square \alpha) \leq h(\alpha)$, for all \,  $h \in Hom(\mathfrak{Fm}, \mathfrak{M}_{4m})$. But, considering all the cases, we must have that $h(\neg \square \alpha)=0$ \, and \, $h(\alpha)=1$, for all $h$, and therefore the sequent $\Rightarrow \alpha$ is provable in $\mathfrak{G}$. The converse is straightforward.
\end{dem}

\

\noindent Recall that a rule of inference is {\em admissible} in a
formal system if the set of theorems of the system is closed
under the rule; and a rule is said to be {\em derivable} in the
same formal system if its conclusion can be derived from its
premises using the other rules of the system.

\noindent A well--known rule for readers familiar with modal logic is the {\em Rule of Necessitation},
which states that if $\varphi$ is a theorem, so is $\square \varphi$. Formally,

$$ \mbox{(Nec) \,} \frac{ \Rightarrow \varphi } { \Rightarrow \square \varphi} $$

\noindent Then, we have that:

\begin{lem}\label{LemNec} The Rule of Necessitation is admissible in $\mathfrak{G}$.
\end{lem}
\begin{dem} From Corollary \ref{CorTheoiffValid} and considering the algebra $\mathfrak{M}_{4m}$.
\end{dem}

\

\noindent From the above lemma, we can obtain a proof of $\Rightarrow\square (\alpha \vee \neg \square \alpha)$ in
$\mathfrak{G}$, for any $\alpha \in Fm$. Let $\Pi$ be
a proof of $\Rightarrow \square (\alpha \vee \neg \square \alpha)$ and
let ($r$) be the last rule application in $\Pi$. Clearly, $\Pi$ make use of
more than one rule since $\square (\alpha \vee \neg \square
\alpha)$ is not an axiom. Then, we have the following two cases:

\begin{prooftree}
\AxiomC{Case 1: $\Pi$ is of the form}
\noLine
\UnaryInfC{$\cdot$}
\noLine
\UnaryInfC{$\cdot$}
\noLine
\UnaryInfC{$\cdot$}
\noLine
\UnaryInfC{$\Gamma \Rightarrow \varphi$}
\LeftLabel{\small(r)}
\UnaryInfC{$\Rightarrow \square (\alpha \vee \neg \square \alpha)$}

\AxiomC{}

\AxiomC{Case 2: $\Pi$ is of the form}
\noLine
\UnaryInfC{$\cdot$}
\noLine
\UnaryInfC{$\cdot$}
\noLine
\UnaryInfC{$\cdot$}
\noLine
\UnaryInfC{$\Gamma_1 \Rightarrow \varphi_1 \, \, \, \, \Gamma_2 \Rightarrow \varphi_2$}
\LeftLabel{\small(r)}
\UnaryInfC{$\Rightarrow \square (\alpha \vee \neg \square \alpha)$}
\noLine
\TrinaryInfC{}
\end{prooftree}

\noindent In case 1, ($r$) has just one premise, and therefore it
can be: ($\bot$), weakening,  ($\wedge \Rightarrow$), ($\vee \Rightarrow$),
 ($\Rightarrow \vee$), ($\neg$), ($\neg \neg \Rightarrow$), ($\Rightarrow \neg \neg $),  ($\square \Rightarrow$) or ($\Rightarrow \square$). In the case of ($\bot$), the only possibility is having $\Gamma=\emptyset$. But this would imply that the sequent $\Rightarrow \bot$ is provable, which contradicts the soundness of $\mathfrak{G}$. Thus, this case is discarded. On the other hand, none of the other rules above has the structure of ($r$), so they are also discarded.\\
Therefore, $\pi$ is of the form depicted in Case 2. Then, ($r$) must be one of the following: the cut rule, ($\Rightarrow \wedge$) or ($\vee \Rightarrow$).
It is clear that ($r$) cannot be ($\Rightarrow \wedge$) nor ($\vee \Rightarrow$). Consequently, ($r$) must be the cut rule.

\

\noindent We have just proved, therefore, the following assertion.

\begin{prop} Every proof of \, $\Rightarrow \square (\alpha \vee \neg \square \alpha)$  in  $\mathfrak{G}$ uses the cut rule.
\end{prop}

\noindent Moreover, we have that:

\begin{lem} For every $\varphi \in Fm$ such that $\Rightarrow \varphi$  is provable in $\mathfrak{G}$, we have that $\Rightarrow \square \varphi$ is provable in $\mathfrak{G}$; and every proof of $\Rightarrow \square \varphi$ in $\mathfrak{G}$ makes use of the cut rule.
\end{lem}

\noindent Consequently,

\begin{theo} $\mathfrak{G}$ does not admit cut--elimination.
\end{theo}

\section{The general method of Avron, Ben-Naim and Konikowska} \label{s2}

In \cite{Avron01}, A. Avron and B. Konikowska use the Rasiowa-Sikorski  decomposition
methodology to get sound and complete proof systems employing $n$-sequents
for all propositional logics based on non-deterministic matrices. Later, these same authors jointly with J. Ben-Naim (\cite{Avron02}) presented a general method to transform a 
given sound and complete $n$-sequent proof system into an equivalent sound and complete system of ordinary two-sided sequents (for languages satisfying a certain minimal expressiveness condition). In this section  we shall recall both methods considering ordinary (deterministic) matrices.

In what follows, $\mathscr{L}$ is a propositional language and let (in this section) $\mathfrak{Fm}$ be the absolutely free algebra over $\mathscr{L}$ generated by some denumerable set of variables, with underlying set (of formulas) $Fm$. Let ${\cal M}=\langle {\cal T}, {\cal D }, {\cal O}\rangle$ be a logical matrix for $\mathscr{L}$. As we said, a valuation $v$ in ${\cal M}$ satisfies a given formula $\alpha$ if $v(\alpha)\in {\cal D}$.
A sequent $\Gamma \Rightarrow \Delta$ is satisfied by the valuation $v$, in symbols $v\models \, \Gamma \Rightarrow \Delta$, if either $v$ does not satisfy some formula in $\Gamma$ or $v$ satisfies some formula in $\Delta$. A sequent is {\em valid} if it is satisfied by all valuations. 
 
Now, suppose that ${\cal T}=\{t_0,\dots, t_{n-1}\}$, where $n\geq 2$,  and ${\cal D}=\{t_d,\dots, t_{n-1}\}$, where $1\leq d \leq n-1$.

\begin{defi} (see \cite{Avron01}) An $n$--sequent over $\mathscr{L}$ is an expression  
$$ \Gamma_0 \mid \dots \mid \Gamma_{n-1}$$
where, for each $i$, $\Gamma_i$ is a finite set of formulas. A valuation $v$ satisfies the $n$--sequent $\Gamma_0 \mid \dots \mid \Gamma_{n-1}$ iff there exists $i$, $0\leq i \leq n-1$ \, and $\psi\in \Gamma_i$ such that $v(\psi)=t_i$. An $n$--sequent is valid if it is satisfied by every valuation $v$. 
\end{defi}  

\noindent Note that, a valuation $v$ satisfies an ordinary sequent $\Gamma  \Rightarrow \Delta$ iff $v$ satisfies the $n$--sequent $\Gamma_1 \mid \dots \mid \Gamma_{n-1}$ where $\Gamma_i=\Gamma$ for all $0\leq i\leq d-1$ and $\Gamma_j=\Delta$ for all $d\leq j\leq n-1$ .

\

An alternative presentation of $n$-sequents is by means of sets of signed formulas. A signed formula over the language $\mathscr{L}$ and ${\cal T}$, is an expression of the form
$$t_i : \psi$$
where $t_i\in  {\cal T}$ and $\psi \in Fm$. A valuation $v$ satisfies the signed formula $t_i : \psi$ iff $v(\psi)=t_i$. If ${\cal U}\subseteq {\cal T}$ and $\Gamma \subseteq Fm$, we denote by ${\cal U} : \Gamma$ the set

$${\cal U} : \Gamma =\{ t:\alpha \mid t\in {\cal U}, \alpha \in \Gamma\}$$

If ${\cal U}=\{t\}$, we write $t : \Gamma$ instead of $\{t\} : \Gamma$.  
A valuation satisfies the set of signed formulas ${\cal U} : \Gamma$  if it satisfies some signed formula of ${\cal U} : \Gamma$; and we say that ${\cal U} : \Gamma$ is valid if it is satisfied by every valuation $v\in {\cal V}$.
It is clear that,  the $n$--sequent $ \Gamma_0 \mid \dots \mid \Gamma_{n-1}$ is valid iff the set of signed formulas $\bigcup \limits_{i=0}^ {n-1} t_i : \Gamma_i$ is valid.

\

A. Avron and B. Konikowska developed in \cite{Avron01} a generic $n$-sequent system for any logic based on an $n$-valued matrix. Consider the $n$-valued matrix ${\cal M}=\langle {\cal T}, {\cal D }, {\cal O}\rangle$ and let $SF_{\cal M}$ the system defined as follows: for $\Omega$ and $\Omega'$ sets of signed formulas

\begin{itemize}
\item {\bf Axioms:} \, \,    ${\cal T}: \alpha$
\item {\bf Structural rules: } Weakening:

$$ \displaystyle  \frac{\Omega } {\Omega'} \hspace{0.3cm} \mbox{ in case } \hspace{0.2cm} \Omega\subseteq \Omega'$$

\item {\bf Logical rules: } for each $k$-ary connective $f$ and every $(a_1, \dots, a_k) \in   {\cal T}^k$ 
$$ \displaystyle  \frac{\Omega, a_1: \alpha_1 \, \, \dots \, \, \Omega, a_k: \alpha_k } {\Omega,  \hat{f}(a_1, \dots, a_k) : f(\alpha_1, \dots,\alpha_k)}$$
\end{itemize}

\

\begin{theo}\label{TeoAv01} (\cite{Avron01}) The system $SF_{\cal M}$ is sound and complete w.r.t. the matrix ${\cal M}$
\end{theo}

Let $Fm_p$ be the set of all formulas of  $Fm$ that have $p$ as their only propositional variable, i.e., $Fm_p=\{ \alpha  \in Fm : Var(\alpha)=\{p\} \}$. Let ${\cal M}=\langle {\cal T}, {\cal D }, {\cal O}\rangle$ be a logical matrix and denote by ${\cal N}$ the set ${\cal T}\setminus {\cal D}$.

\begin{defi}(\cite{Avron02}) \label{SufficientlyExp} The language $\mathscr{L}$ is sufficiently expressive for ${\cal M}$ iff for any $i$, $0\leq i \leq n-1$ there exist natural numbers $l_i,m_i$  and formulas $\alpha_{j}^{i}, \beta_{k}^{i}\in Fm_p$, for $1\leq j\leq l_i$ and $1\leq k\leq m_i$ such that for any valuation $v$, the following conditions hold:\\
(i) $\alpha_{1}^{i}=p$ if $t_i\in {\cal N}$ \, and \, $\beta_{1}^{i}=p$ if $t_i\in {\cal D}$, \\
(ii) For $\varphi \in Fm$ and $t_i\in {\cal T}$
$$v(\varphi)=t_i \, \Leftrightarrow \ v(\alpha_{1}^{i}[p/\varphi]), \dots, v(\alpha_{l_i}^{i}[p/\varphi]) \in {\cal N} \, and \,  v(\beta_{1}^{i}[p/\varphi]), \dots, v(\alpha_{m_i}^{i}[p/\varphi]) \in {\cal D}$$
where $\alpha_{j}^{i}[p/\varphi]$ ($\beta_{k}^{i}[p/\varphi]$) is the formula obtained by the substitution of $p$ by $\varphi$ in $\alpha_{j}^{i}$ ($\beta_{k}^{i}$).
\end{defi}

\noindent Note that, as it is mentioned in \cite{Avron02}, condition (i) above is not really limiting, since given $\alpha_{j}^{i}, \beta_{k}^{i}$
satisfying (ii), we can simply add to them the necessary formula $p$ without violating (ii). Condition (i) will only be used for a backward translation from ordinary sequents to $n$-sequents, and will be disregarded otherwise.

\

\noindent If $\Gamma$ is a set of formulas and $\alpha\in Fm_p$, we denote by $\alpha[p/\Gamma]$ the set

$$\alpha[\Gamma]=\{\alpha[p/\gamma] \mid \gamma\in \Gamma\}$$
The method is based on replacing each $n$-sequent by a semantically equivalent set of two-sided sequents.

\

\noindent Let $\mathscr{L}$ be a sufficiently expressive language and for $0\leq i \leq n-1$ let $l_i$, $m_i$, $\alpha_{j}^{i}$  and $\beta_{k}^{i}$ as in Definition \ref{SufficientlyExp}. Consider the $n$--sequent $\Sigma = \Gamma_0 \mid \dots \mid \Gamma_{n-1}$ over $\mathscr{L}$. A partition $\pi$ of the $n$--sequent $\Sigma$ is a tuple $\pi=(\pi_0,\dots,\pi_{n-1})$  such that, for every $i$, $\pi_i$ is a partition of the set $ \Gamma_i$ of the form:

$$\pi_i=\{\Gamma'_{ij}\mid 1\leq j\leq l_i\}\cup\{\Gamma''_{ik}\mid 1\leq k\leq m_i\}$$

\noindent  Note that $\pi_i$ is not a partition in the usual sense, since its components are allowed to be empty. Besides, observe that the number of sets in this partition is exactly the number of formulas corresponding to $i$ in Definition \ref{SufficientlyExp}.

\

Then, given a partition $\pi$ of the $n$-sequent $\Sigma$, we define the two-sided sequent $\Sigma_{\pi}$ determined by $\Sigma$ and the partition $\pi$,  as follows:

$$\bigcup \limits_{j=0}^{l_0} \alpha_{j}^{0}[\Gamma'_{0j}], \dots, \bigcup \limits_{j=n-1}^{l_{n-1}} \alpha_{j}^{n-1}[\Gamma'_{(n-1)j}] \, \Rightarrow \, \bigcup \limits_{k=0}^{m_0} \beta_{k}^{0}[\Gamma''_{0k}], \dots, \bigcup \limits_{k=n-1}^{m_{n-1}} \beta_{k}^{n-1}[\Gamma''_{(n-1)k}]  $$

\noindent Let $\Pi$ be the set of all partitions of the $n$--sequent $\Sigma$. Then, the set $TWO(\Sigma)$ is defined as follows: 

$$TWO(\Sigma)= \{\Sigma_{\pi}\mid \pi\in \Pi\}$$

\begin{theo}(\cite{Avron02}) Let $\Sigma$ be an $n$--sequent over $\mathscr{L}$ and $v$ a valuation. Then,
$v$ satisfies $\Sigma$ \,iff \, $v$ satisfies $\Sigma'$, for every $\Sigma' \in TWO(\Sigma)$. 
\end{theo}

\begin{defi}(\cite{Avron02}) \label{defGC} Let ${\cal C}$ be an  $n$--sequent calculus over $\mathscr{L}$. Then, let  $TWO({\cal C})$ the (ordinary) sequent calculus over $\mathscr{L}$ given by:
\begin{itemize}
\item[] {\bf Axioms:} $TWO(A)$, for all axiom $A$ of $\cal C$,
\item[] {\bf Inference rules:} $\displaystyle \frac{TWO(S)}{\Sigma'}$, where $S$ is a finite set of $n$-sequents, $R$ is one $n$-sequent such that $\displaystyle \frac{S}{R}$ is a rule in $\cal C$ and $\Sigma'\in TWO(R)$.
\end{itemize}
\end{defi}

\

\noindent Then, 
 
\begin{theo}(\cite{Avron02})\label{teoSeq01} If an $n$--sequent $\Sigma$ is provable in $\cal C$, then each two-sided sequent $\Sigma'\in TWO(\Sigma)$ is provable in $TWO({\cal C})$.
\end{theo}

\begin{theo}\label{TeoAv3.5}(\cite{Avron02})\label{teoTWO} Let $\mathscr{L}$ be a sufficiently expressive language for ${\cal M}$, and let $\cal C$ be a sound and complete sequent calculus w.r.t $\cal M$. Then, $TWO({\cal C})$ is sound and complete w.r.t. ${\cal M}$.
\end{theo}

\noindent The analogue of the cut rule for ordinary sequents is the following generalized cut rule for sets of signed formulas:

$$\displaystyle \frac{\Omega \cup \{i:\alpha \, | \, i\in I\} \hspace{0.5cm}\Omega \cup \{j:\alpha \, | \, j\in J\}}{\Omega}\hspace{0.5cm} \mbox{ for } I, J\subseteq {\cal V}, I\cap J=\emptyset$$

\begin{theo}(\cite{Avron02}) Under the conditions of Theorem \ref{teoTWO}, the cut rule is admissible in $TWO({\cal C})$. In particular, if $\cal C$ is obtained by the method of \cite{Avron01}, then the cut rule is admissible in $TWO({\cal C})$.
\end{theo}

As it was observed in \cite{Avron02}, the $n$-sequent calculi obtained using the above general method are hardly optimal (the same is true for the two-sided calculi). We can use the three general streamlining principles from \cite{Avron01} to reduce the calculi to a more compact form.  The three streamlining principles are: {\em Principle 1:} deleting a derivable rule, {\em Principle 2:} simplifying a rule by replacing it with one with weaker premises, and {\em Principle 3:} combining two context--free rules with the same conclusion into one. Recall that a rule $R$ is context-free if whenever $\frac{\phi_1\dots\phi_n}{\Sigma}$ is a valid application of $R$, and $\Sigma'$ is a set of signed formulas, then $\frac{\phi_1\cup\Sigma'\dots\phi_n\cup\Sigma'}{\Sigma\cup\Sigma'}$ is also a valid application of $R$. A rule $R$ of an ordinary two--sided sequent calculus is a context--free if  \, $\displaystyle \frac{\Gamma_1\Rightarrow\Delta_1, \dots, \Gamma_k\Rightarrow\Delta_k}{\Gamma\Rightarrow\Delta}$ \, is a valid application of $R$, then \,   $\displaystyle \frac{\Gamma_1, \Gamma' \Rightarrow\Delta_1, \Delta', \dots, \Gamma_k,\Gamma'\Rightarrow\Delta_k, \Delta'}{\Gamma,\Gamma'\Rightarrow\Delta, \Delta'}$ \, is also a valid application of $R$, where $\Gamma'$ and $\Delta'$ are finite sets of formulas. \\
Of these three, the first and the third decrease the number of rules, while the second simplifies a rule by decreasing the number of its premises.

It is worth mentioning that applying Principles 1–3 preserves the cut-elimination property
since cut-elimination is obtained via the completeness result and the principles are designed to
retain completeness.

\

\section{Cut--free sequent calculus for ${\cal TML}$}

Now, we shall use the method exhibited in Section \ref{s2} to develop a $4$-sequent calculus for $\cal TML$.  In this case, we shall use its alternative presentation provided by sets of $4$-signed formulas. 

\noindent Let ${\cal SF}_4$ be $4$-sequent calculus given by: for $\alpha, \beta \in Fm$, $\Omega$ and $\Omega'$ arbitrary sets of signed formulas

\

\noindent{\bf Axioms:} \, \, \,  $\{0:\alpha, {\bf n}:\alpha, {\bf b}:\alpha,  n:\alpha \}$.

\

\noindent{\bf Structural rules:} Weakening.

$$ \displaystyle  \frac{\Omega } {\Omega'} \hspace{0.3cm} \mbox{ in case } \hspace{0.2cm} \Omega\subseteq \Omega'$$
\

\noindent{\bf Logical rules:}  for $i, j \in M_4$

$$ \mbox{($\vee_{ij}$) \, } \displaystyle \frac{\Omega, i: \alpha \hspace{1cm} \Omega, j: \beta } {\Omega, \Sup\{i,j\}: \alpha \vee \beta} \hspace{2.5cm} \mbox{($\wedge_{ij}$) \, } \displaystyle \frac{\Omega, i: \alpha \hspace{1cm} \Omega, j: \beta } {\Omega, \Inf\{i,j\}: \alpha \wedge \beta} $$

$$ \mbox{($\neg_0$) \, } \displaystyle \frac{\Omega, 0: \alpha } {\Omega, 1: \neg \alpha} \hspace{1cm} \mbox{($\neg_{\bf n}$) \, } \displaystyle \frac{\Omega, {\bf n}: \alpha } {\Omega, {\bf n}: \neg \alpha} \hspace{1cm} 
\mbox{($\neg_{\bf b}$) \, } \displaystyle \frac{\Omega, {\bf b}: \alpha } {\Omega, {\bf b}: \neg \alpha} \hspace{1cm}
\mbox{($\neg_1$) \, } \displaystyle \frac{\Omega, 1: \alpha } {\Omega, 0: \neg \alpha}$$

$$ \mbox{($\square_{i}$) \, } \displaystyle \frac{\Omega, i: \alpha } {\Omega, 0: \square \alpha},  
 \mbox{ for } i\not=1 \hspace{1.5cm} \mbox{($\square_1$) \, } \displaystyle \frac{\Omega, 1: \alpha } {\Omega, 1: \square \alpha} $$

\

\

\noindent In rules ($\vee_{ij}$) (and (($\wedge_{ij}$)), the supremum (infimum) is taken on the lattice $M_4$. Besides, observe that the system ${\cal SF}_4$ has forty logical rules and it is not optimal. However, in this step we are not going to use the principles mentioned in Section \ref{s2} to reduce ${\cal SF}_4$.

\begin{prop}\begin{itemize} \item[]
\item[(i)] ${\cal SF}_4$ is sound and complete w.r.t. the matrix ${\cal M}_4$,
\item[(ii)] the cut rule is admissible in ${\cal SF}_4$.
\end{itemize}
\end{prop}     
\begin{dem} From Theorem \ref{TeoAv01}.
\end{dem}

\

\noindent Now, we shall apply the method described in Section \ref{s2} to translate ${\cal SF}_4$  to an ordinary two-sided sequent calculus.

\begin{prop} The language $\mathscr{L}$ is sufficiently expressive for the semantics determined by the matrix ${\cal M}_4$.
\end{prop}     
\begin{dem} Let $v:Fm \to M_4$ be a valuation and let $\alpha\in Fm$ an arbitrary formula, then we have  that

$$v(\alpha)=0  \, \Longleftrightarrow \, v(\alpha)\in {\cal N} \mbox{ and } v(\neg \alpha)\in {\cal D}$$  
$$v(\alpha)={\bf n} \, \Longleftrightarrow \, v(\alpha)\in {\cal N} \mbox{ and } v(\neg \alpha)\in {\cal N}$$ 
$$v(\alpha)={\bf b}  \, \Longleftrightarrow \, v(\alpha)\in {\cal D} \mbox{ and } v(\neg \alpha)\in {\cal D}$$
$$v(\alpha)=1  \, \Longleftrightarrow \, v(\alpha)\in {\cal D} \mbox{ and } v(\neg \alpha)\in {\cal N}$$     
where ${\cal N}=M_4\setminus {\cal D}=\{0, {\bf n}\}$. 
\end{dem}

\

\noindent  According to Theorem \ref{TeoAv3.5}, to transform ${\cal SF}_4$ to an ordinary one, we have to replace every axiom $A$ with the equivalent set of ordinary sequents $TWO(A)$. In terms of $4$-sequents, the only axiom of ${\cal SF}_4$ has the form

$$\alpha\mid \alpha\mid \alpha\mid \alpha$$
and it yields the following ordinary two-sided sequents
$$\alpha, \neg\alpha\Rightarrow \alpha \hspace{1cm} \alpha, \neg\alpha\Rightarrow \neg\alpha \hspace{1cm} \alpha\Rightarrow \neg\alpha, \alpha \hspace{1cm} \alpha, \neg\alpha\Rightarrow \neg\alpha, \alpha \hspace{1cm} \neg\alpha\Rightarrow \neg\alpha, \alpha $$
All of them can be derived from \, $\alpha \Rightarrow \alpha$ \, (or from an instance of it) by the use of weakening. \\
Now, let us focus on  rules  ($\vee_{ij}$), $i,j\in M_4$. First observe that, if $\varphi\in Fm$ then

$$TWO(\varphi\mid \hspace{.3cm} \mid \hspace{.3cm} \mid \hspace{.3cm} )=\{\varphi\Rightarrow \, , \Rightarrow \neg \varphi\}$$   
$$TWO(\hspace{.3cm} \mid \varphi \mid \hspace{.3cm} \mid \hspace{.3cm} )=\{\varphi\Rightarrow \, ,  \neg \varphi\Rightarrow\}$$   
$$TWO(\hspace{.3cm} \mid \hspace{.3cm} \mid \varphi \mid \hspace{.3cm} )=\{\Rightarrow \varphi \, ,  \Rightarrow \neg \varphi\}$$   
$$TWO(\hspace{.3cm} \mid \hspace{.3cm} \mid \hspace{.3cm} \mid \varphi )=\{ \neg \varphi \Rightarrow\, ,  \Rightarrow \varphi\}$$

\noindent So, after removing the contexts for brevity, the rules ($\vee_{ij}$)'s are translated to the following thirty-two two-sided sequent rules: 

\

\begin{tabular}{llll}
($\vee$)$_{10}$ & $\displaystyle \frac{\Rightarrow \alpha \hspace{.3cm} \neg\alpha\Rightarrow \hspace{.3cm} \beta\Rightarrow \hspace{.3cm} \Rightarrow \neg\beta}{\Rightarrow \alpha \vee\beta \mbox{ \, ; \, } \neg(\alpha \vee \beta) \Rightarrow}$ & ($\vee$)$_{1{\bf n}}$ & $\displaystyle \frac{\Rightarrow \alpha \hspace{.3cm} \neg\alpha\Rightarrow \hspace{.3cm} \beta\Rightarrow \hspace{.3cm} \neg\beta \Rightarrow}{\Rightarrow \alpha \vee\beta \mbox{ \, ; \, } \neg(\alpha \vee \beta) \Rightarrow}$ \\

\\		

($\vee$)$_{1{\bf b}}$ & $\displaystyle \frac{\Rightarrow \alpha \hspace{.3cm} \neg\alpha\Rightarrow \hspace{.3cm} \Rightarrow\beta \hspace{.3cm} \Rightarrow \neg\beta}{\Rightarrow \alpha \vee\beta \mbox{ \, ; \, } \neg(\alpha \vee \beta) \Rightarrow}$ & ($\vee$)$_{11}$ & $\displaystyle \frac{\Rightarrow \alpha \hspace{.3cm} \neg\alpha\Rightarrow \hspace{.3cm} \Rightarrow\beta \hspace{.3cm} \neg\beta \Rightarrow}{\Rightarrow \alpha \vee\beta \mbox{ \, ; \, } \neg(\alpha \vee \beta) \Rightarrow}$ \\		

\\

($\vee$)$_{{\bf b}0}$ & $\displaystyle \frac{\Rightarrow \alpha \hspace{.3cm} \neg\alpha\Rightarrow \hspace{.3cm} \beta\Rightarrow \hspace{.3cm} \Rightarrow \neg\beta}{\Rightarrow \alpha \vee\beta \mbox{ \, ; \, } \Rightarrow\neg(\alpha\vee\beta)}$ & ($\vee$)$_{{\bf b}{\bf n}}$ & $\displaystyle \frac{\Rightarrow \alpha \hspace{.3cm} \neg\alpha\Rightarrow \hspace{.3cm} \beta\Rightarrow \hspace{.3cm} \neg\beta \Rightarrow}{\Rightarrow \alpha \vee\beta \mbox{ \, ; \, } \neg(\alpha \vee \beta)\Rightarrow}$ \\

\\

($\vee$)$_{{\bf b}{\bf b}}$ & $\displaystyle \frac{\Rightarrow \alpha \hspace{.3cm} \neg\alpha\Rightarrow \hspace{.3cm} \Rightarrow\beta \hspace{.3cm} \Rightarrow \neg\beta}{\Rightarrow \alpha \vee\beta \mbox{ \, ; \, } \Rightarrow\neg(\alpha\vee\beta)}$ & ($\vee$)$_{{\bf b}1}$ & $\displaystyle \frac{\Rightarrow \alpha \hspace{.3cm} \neg\alpha\Rightarrow \hspace{.3cm} \Rightarrow\beta \hspace{.3cm} \neg\beta \Rightarrow}{\Rightarrow \alpha \vee\beta \mbox{ \, ; \, } \neg(\alpha \vee \beta)\Rightarrow}$ \\		
\end{tabular}

\begin{tabular}{llll}
($\vee$)$_{{\bf n}0}$ & $\displaystyle \frac{\alpha\Rightarrow \hspace{.3cm} \neg\alpha\Rightarrow \hspace{.3cm} \beta\Rightarrow \hspace{.3cm} \Rightarrow \neg\beta}{\alpha \vee\beta \Rightarrow \mbox{ \, ; \, } \neg(\alpha \vee \beta) \Rightarrow}$ & ($\vee$)$_{{\bf n}{\bf n}}$ & $\displaystyle \frac{\alpha\Rightarrow \hspace{.3cm} \neg\alpha\Rightarrow \hspace{.3cm} \beta\Rightarrow \hspace{.3cm} \neg\beta \Rightarrow}{ \alpha \vee\beta \Rightarrow \mbox{ \, ; \, } \neg(\alpha \vee \beta) \Rightarrow}$ \\

\\		

($\vee$)$_{{\bf n}{\bf b}}$ & $\displaystyle \frac{\alpha\Rightarrow \hspace{.3cm} \neg\alpha\Rightarrow \hspace{.3cm} \Rightarrow\beta \hspace{.3cm} \Rightarrow \neg\beta}{\Rightarrow \alpha \vee\beta \mbox{ \, ; \, } \neg(\alpha \vee \beta) \Rightarrow}$ & ($\vee$)$_{{\bf n}1}$ & $\displaystyle \frac{\alpha\Rightarrow \hspace{.3cm} \neg\alpha\Rightarrow \hspace{.3cm} \Rightarrow\beta \hspace{.3cm} \neg\beta \Rightarrow}{\Rightarrow \alpha \vee\beta \mbox{ \, ; \, } \neg(\alpha \vee \beta) \Rightarrow}$ \\		

\\

($\vee$)$_{00}$ & $\displaystyle \frac{\alpha\Rightarrow \hspace{.3cm} \neg\alpha\Rightarrow \hspace{.3cm} \beta\Rightarrow \hspace{.3cm} \Rightarrow \neg\beta}{\alpha \vee\beta \Rightarrow \mbox{ \, ; \, } \Rightarrow\neg(\alpha\vee\beta)}$ & ($\vee$)$_{0{\bf n}}$ & $\displaystyle \frac{\alpha\Rightarrow \hspace{.3cm} \neg\alpha\Rightarrow \hspace{.3cm} \beta\Rightarrow \hspace{.3cm} \neg\beta \Rightarrow}{\alpha \vee\beta \Rightarrow \mbox{ \, ; \, } \neg(\alpha \vee \beta)\Rightarrow}$ \\

\\		

($\vee$)$_{0{\bf b}}$ & $\displaystyle \frac{\alpha\Rightarrow \hspace{.3cm} \neg\alpha\Rightarrow \hspace{.3cm} \Rightarrow\beta \hspace{.3cm} \Rightarrow \neg\beta}{\Rightarrow \alpha \vee\beta \mbox{ \, ; \, } \Rightarrow\neg(\alpha\vee\beta)}$ & ($\vee$)$_{01}$ & $\displaystyle \frac{\alpha\Rightarrow \hspace{.3cm} \neg\alpha\Rightarrow \hspace{.3cm} \Rightarrow\beta \hspace{.3cm} \neg\beta \Rightarrow}{\Rightarrow \alpha \vee\beta \mbox{ \, ; \, }  \neg(\alpha \vee \beta)\Rightarrow}$ \\		

\end{tabular}
		
\	

\

\noindent In the above list  we use an informal notation by separating the alternate conclusion sequents with semicolons. At this point, we shall follow the three principles mentioned in the above section in order to reduce the number of rules. Our main tool for this job will be the next proposition. 

\

\begin{prop}\label{propAux} Let $\mathfrak{SC}$ a sequent calculus in which the cut rule is admissible,  let $S$ be a set of sequents and $\Sigma$ be a sequent such that $\displaystyle \frac{S \cup \{\Gamma\Rightarrow \Delta, \varphi\}}{\Sigma}$ and $\displaystyle \frac{S \cup \{\Gamma, \varphi\Rightarrow\Delta\}}{\Sigma}$ are two context-free rules of $\mathfrak{SC}$. Then,  $\displaystyle \frac{S}{\Sigma}$ is derivable in $\mathfrak{SC}$.
\end{prop}
\begin{dem} From the fact that the rules are context-free and using the cut rule.
\end{dem}\\[3mm]
Then, from ($\vee$)$_{10}$, ($\vee$)$_{1{\bf n}}$ and Proposition \ref{propAux} we get  \, $\displaystyle \frac{\Rightarrow \alpha \hspace{.3cm} \neg\alpha\Rightarrow \hspace{.3cm} \beta\Rightarrow}{\Rightarrow \alpha \vee\beta \mid \neg(\alpha \vee \beta) \Rightarrow}$. From ($\vee$)$_{1{\bf b}}$, ($\vee$)$_{11}$ and Proposition \ref{propAux} we get \, $\displaystyle \frac{\Rightarrow \alpha \hspace{.3cm} \neg\alpha\Rightarrow \hspace{.3cm} \Rightarrow\beta}{\Rightarrow \alpha \vee\beta\mid \neg(\alpha \vee \beta) \Rightarrow}$.\\[1.3mm]
 From these rules and Proposition \ref{propAux} we obtain (1) $\displaystyle \frac{\Rightarrow \alpha \hspace{.3cm} \neg\alpha\Rightarrow}{\Rightarrow \alpha \vee\beta}$ \, and \, (1') $\displaystyle \frac{\Rightarrow \alpha \hspace{.3cm} \neg\alpha\Rightarrow}{\neg(\alpha \vee \beta) \Rightarrow}$.
Analogously, from ($\vee$)$_{{\bf b}0}$, ($\vee$)$_{{\bf b}{\bf n}}$, ($\vee$)$_{{\bf b}{\bf b}}$, ($\vee$)$_{{\bf b}1}$ we obtain\\[1.3mm] (2)
$\displaystyle \frac{\Rightarrow \alpha \hspace{.3cm} \Rightarrow\neg\alpha}{\Rightarrow \alpha \vee\beta}$. Finally, from (1), (2) and Proposition \ref{propAux} we get that 
\begin{equation}\tag{3}\label{eq3}
\displaystyle \frac{\Rightarrow \alpha}{\Rightarrow \alpha \vee\beta}
\end{equation}
is derivable. On the other hand, following an analogous reasoning we can prove that 
\begin{equation}\tag{4}\label{eq4}
\displaystyle \frac{\Rightarrow \beta}{\Rightarrow \alpha \vee\beta}
\end{equation}
is derivable. Then, after combining rules (\ref{eq3}) and (\ref{eq4}) and restoring the context we get the rule
$$\mbox{($\Rightarrow\vee$) } \displaystyle \frac{\Gamma \Rightarrow \Delta, \alpha, \beta}{\Gamma \Rightarrow \Delta, \alpha\vee\beta}$$
From ($\vee$)$_{{\bf n}0}$, ($\vee$)$_{{\bf n}{\bf n}}$, ($\vee$)$_{{\bf n}{\bf b}}$ and ($\vee$)$_{{\bf n}1}$ we obtain (5) \, $\displaystyle \frac{\Rightarrow \alpha \hspace{.3cm} \neg\alpha\Rightarrow}{\Rightarrow \alpha \vee\beta}$; then using (1') and restoring the context we get (5) $\displaystyle \frac{\Gamma, \neg\alpha\Rightarrow \Delta}{\Gamma, \neg(\alpha \vee \beta) \Rightarrow \Delta}$. In a similar way, it can be proved that  (6) $\displaystyle \frac{\Gamma, \neg\beta\Rightarrow \Delta}{\Gamma, \neg(\alpha \vee \beta) \Rightarrow \Delta}$ is derivable. Then, combining (5) and (6) and restoring the context we get
$$\mbox{($\neg\vee\Rightarrow$) }\displaystyle \frac{\Gamma, \neg\alpha, \neg\beta\Rightarrow \Delta}{\Gamma, \neg(\alpha \vee \beta) \Rightarrow \Delta}$$
From ($\vee$)$_{{\bf n}0}$, ($\vee$)$_{{\bf n}{\bf n}}$, ($\vee$)$_{00}$ and ($\vee$)$_{0{\bf n}}$ and restoring context we obtain the rule
$$\mbox{($\vee\Rightarrow$) }\displaystyle \frac{\Gamma, \alpha\Rightarrow \Delta \hspace{.5cm }\Gamma,\beta\Rightarrow \Delta}{\Gamma, \alpha \vee \beta \Rightarrow \Delta}$$
and,  from ($\vee$)$_{00}$, ($\vee$)$_{0{\bf b}}$, ($\vee$)$_{{\bf b}0}$ and ($\vee$)$_{{\bf b}{\bf b}}$ we get
$$\mbox{($\Rightarrow\neg\vee$) }\displaystyle \frac{\Gamma\Rightarrow \Delta, \neg\alpha \hspace{.5cm }\Gamma\Rightarrow \Delta, \neg \beta}{\Gamma \Rightarrow \Delta, \neg(\alpha \vee \beta)}$$

\

\noindent In the same way, we obtain the following rules for the connective $\wedge$:

$$ \mbox{($\wedge \Rightarrow$) \, } \displaystyle \frac{\Gamma, \alpha, \beta \Rightarrow \Delta} {\Gamma, \alpha \wedge \beta \Rightarrow \Delta} \hspace{2cm} \mbox{($\Rightarrow \wedge$) \, }
 \displaystyle \frac{\Gamma \Rightarrow \Delta, \alpha  \hspace{0.5cm} \Gamma \Rightarrow \Delta, \beta}{\Gamma \Rightarrow \Delta, \alpha \wedge \beta} $$

$$ \mbox{($\neg\wedge\Rightarrow$) \, } \displaystyle \frac{\Gamma, \neg\alpha \Rightarrow \Delta  \hspace{0.5cm} \Gamma, \neg\beta \Rightarrow \Delta}{\Gamma, \neg(\alpha \wedge \beta) \Rightarrow \Delta} \hspace{1.5cm} \mbox{($\Rightarrow\neg\wedge $) \, } \displaystyle \frac{\Gamma \Rightarrow \Delta, \neg\alpha, \neg\beta} {\Gamma \Rightarrow \Delta, \neg(\alpha \wedge \beta)}  $$

\

\noindent On the other hand, rules ($\neg$)$_i$ with $i\in M_4$ are translated to (after eliminating the trivial rules)

\

\begin{tabular}{llll}
\hspace{1cm} ($\neg$)$_{0}$ & $\displaystyle \frac{ \alpha\Rightarrow \hspace{.3cm} \Rightarrow\neg\alpha}{ \neg \neg \alpha\Rightarrow}$ \hspace{1cm} & ($\neg$)$_{{\bf n}}$ & $\displaystyle \frac{ \alpha\Rightarrow \hspace{.3cm} \neg\alpha\Rightarrow}{\neg \neg \alpha\Rightarrow}$ \\

\\		

\hspace{1cm} ($\neg$)$_{\bf b}$ & $\displaystyle \frac{ \Rightarrow\alpha \hspace{.3cm} \Rightarrow\neg\alpha}{ \Rightarrow \neg \neg \alpha}$ & ($\neg$)$_{1}$ & $\displaystyle \frac{ \Rightarrow \alpha \hspace{.3cm} \neg\alpha\Rightarrow}{\Rightarrow \neg \neg \alpha}$ \\
\end{tabular}

\

\noindent From ($\neg$)$_{0}$, ($\neg$)$_{{\bf n}}$ and Proposition \ref{propAux} on the one hand; and  ($\neg$)$_{{\bf b}}$, ($\neg$)$_{1}$ and Proposition \ref{propAux} on the other, we obtain

$$ \mbox{($\neg\neg\Rightarrow$) \, } \displaystyle \frac{\Gamma, \alpha \Rightarrow \Delta}{\Gamma, \neg\neg\alpha \Rightarrow \Delta} \hspace{1.5cm} \mbox{($\Rightarrow\neg\neg $) \, } \displaystyle \frac{\Gamma \Rightarrow \Delta, \alpha} {\Gamma \Rightarrow \Delta, \neg\neg\alpha}  $$

\

\noindent Finally, rules ($\square$)$_{i}$ are translated to

\

\begin{tabular}{llll}
\hspace{1cm} ($\square$)$_{0}$ & $\displaystyle \frac{ \alpha\Rightarrow \hspace{.3cm} \Rightarrow\neg\alpha}{ \square\alpha\Rightarrow \mbox{ \, ; \, } \Rightarrow \neg\square\alpha}$ \hspace{1cm} & ($\square$)$_{{\bf n}}$ & $\displaystyle \frac{ \alpha\Rightarrow \hspace{.3cm} \neg\alpha\Rightarrow}{ \square\alpha\Rightarrow \mbox{ \, ; \, } \Rightarrow \neg\square\alpha}$ \\

\\		

\hspace{1cm} ($\square$)$_{\bf b}$ & $\displaystyle \frac{ \Rightarrow\alpha \hspace{.3cm} \Rightarrow\neg\alpha}{ \square\alpha\Rightarrow \mbox{ \, ; \, } \Rightarrow \neg\square\alpha}$ & ($\square$)$_{1}$ & $\displaystyle \frac{ \Rightarrow \alpha \hspace{.3cm} \neg\alpha\Rightarrow}{\Rightarrow \square\alpha  \mbox{ \, ; \, } \neg\square\alpha\Rightarrow}$ \\
\end{tabular}

\

\noindent and, from these rules and Proposition  \ref{propAux}, we obtain

\

$$ \mbox{($\square\Rightarrow$)$_1$ \, } \displaystyle \frac{\Gamma, \alpha \Rightarrow \Delta }{\Gamma, \square\alpha \Rightarrow \Delta} \hspace{1.4cm} \mbox{($\square\Rightarrow$)$_1$ \, } \displaystyle \frac{\Gamma \Rightarrow \Delta, \neg\alpha }{\Gamma, \square\alpha \Rightarrow \Delta} \hspace{1.4cm} \displaystyle  \mbox{($\Rightarrow\square$) \, } \frac{\Gamma \Rightarrow \Delta, \alpha \hspace{.3cm} \Gamma, \neg\alpha \Rightarrow \Delta}{\Gamma \Rightarrow \Delta, \square\alpha} $$
$$  \mbox{($\neg\square\Rightarrow$) \, } \displaystyle \frac{\Gamma \Rightarrow \Delta, \alpha \hspace{.3cm} \Gamma, \neg\alpha \Rightarrow \Delta}{\Gamma, \neg\square\alpha\Rightarrow \Delta}  \hspace{1cm} \mbox{($\Rightarrow\neg\square$)$_1$ \, } \displaystyle \frac{\Gamma, \alpha \Rightarrow \Delta}{\Gamma\Rightarrow \Delta, \neg\square\alpha}   \hspace{.9cm} \mbox{($\Rightarrow\neg\square$)$_2$ \, } \displaystyle \frac{\Gamma \Rightarrow \Delta, \neg\alpha }{\Gamma\Rightarrow \Delta, \neg\square\alpha}$$

\begin{defi} Let ${\bf SC}_{\cal TML}$ be the sequent calculus given by the axiom \, $\alpha\Rightarrow\alpha$ \, the structural rules of cut and left and right weakening 

$$ \mbox{\rm ($w\Rightarrow$) \, } \displaystyle \frac{\Gamma \Rightarrow \Delta}{\Gamma, \alpha \Rightarrow \Delta} \hspace{1.5cm} \mbox{\rm ($\Rightarrow w$) \, } \displaystyle \frac{\Gamma \Rightarrow \Delta} {\Gamma \Rightarrow \Delta, \alpha} $$

\noindent and the logical rules  {\rm ($\vee\Rightarrow$), ($\Rightarrow\vee$), ($\neg\vee\Rightarrow$), ($\Rightarrow\neg\vee$),  ($\wedge\Rightarrow$), ($\Rightarrow\wedge$), ($\neg\wedge\Rightarrow$), ($\Rightarrow\neg\wedge$),  ($\neg\neg\Rightarrow$), ($\Rightarrow\neg\neg$), ($\square\Rightarrow$)$_i$, ($\Rightarrow\square$), ($\neg\square\Rightarrow$), ($\Rightarrow\neg\square$)$_i$ \, $i=1,2$}.
\end{defi}

\

\noindent We shall write $\Gamma\Leftrightarrow\Delta$ to indicate that both the sequents $\Gamma\Rightarrow\Delta$ and $\Delta\Rightarrow\Gamma$ are provable. Then, it is not difficult to verify that \, $\alpha\wedge\neg\alpha\Leftrightarrow\alpha\wedge\neg\square\alpha$, for every formula $\alpha$.  Besides, the modal axiom $\Rightarrow \alpha \vee \neg \square \alpha$ of $\mathfrak{G}$ is derivable in ${\bf SC}_{\cal TML}$. Indeed,

\begin{prooftree}
\AxiomC{$\alpha \Rightarrow \alpha$}
\LeftLabel{\small($\Rightarrow\neg\square$)$_1$}
\UnaryInfC{$\Rightarrow \alpha, \neg\square \alpha$}
\LeftLabel{\small($\Rightarrow\vee$)}
\UnaryInfC{$\Rightarrow \alpha\vee \neg\square \alpha$}
\end{prooftree}

\noindent Moreover, the sequent $\Rightarrow \square(\alpha \vee \neg \square \alpha)$ is derivable in ${\bf SC}_{\cal TML}$ without the cut rule:

\begin{prooftree}
\AxiomC{$\alpha \Rightarrow \alpha$}
\LeftLabel{\small($\Rightarrow\neg\square$)$_1$}
\UnaryInfC{$\Rightarrow \alpha, \neg\square \alpha$}
\LeftLabel{\small($\Rightarrow\vee$)}
\UnaryInfC{$\Rightarrow \alpha\vee \neg\square \alpha$}

\AxiomC{$\neg\alpha \Rightarrow \neg\alpha$}
\LeftLabel{\small($\square\Rightarrow$)$_1$}
\UnaryInfC{$\neg\alpha, \square\alpha\Rightarrow$}
\LeftLabel{\small($\neg\neg\Rightarrow$)}
\UnaryInfC{$\neg\alpha, \neg\neg\square\alpha\Rightarrow$}
\LeftLabel{\small($\neg\vee\Rightarrow$)}
\UnaryInfC{$\neg(\alpha\vee\neg\square\alpha)\Rightarrow$}

\LeftLabel{\small($\Rightarrow\square$)}
\BinaryInfC{$\Rightarrow \square(\alpha\vee\neg\square \alpha)$}
\end{prooftree}

\ 

\begin{rem}\label{RemBot}  In Font and Rius' system $\mathfrak{G}$ , the propositional constant $\bot$ is used. By following  Avron, Ben-Naim and Konikowska's method, we obtained a system in which $\bot$ does not appear. However, it is easy to check that  the sequent $\neg \alpha\wedge \square \alpha \Rightarrow$ is provable in ${\bf SC}_{\cal TML}$, for any formula $\alpha$. Then, if we denote by $\bot$ the formula  $\neg \alpha\wedge \square \alpha$, for any formula $\alpha$, we have that the rule $(\bot)$ of   $\mathfrak{G}$ is derivable in ${\bf SC}_{\cal TML}$.
\end{rem}

\begin{theo} \label{TeoComSC}
\begin{itemize}
\item[]
\item[(i)] ${\bf SC}_{\cal TML}$ is sound and complete w.r.t. ${\cal M}_4$.
\item[(ii)] The cut rule is admissible in ${\bf SC}_{\cal TML}$,
\end{itemize}
\end{theo}
\begin{dem} The system ${\bf SC}_{\cal TML}$ was constructed according to the method displayed in Section \ref{s2}.
\end{dem}

\begin{cor}\label{CoroCut} ${\bf SC}_{\cal TML}$ is a cut-free sequent calculus that provides a syntactical counterpart for $\cal TML$. 
\end{cor}

\section{Some applications of the cut elimination theorem}

In this section, we shall use the cut-free system ${\bf SC}_{\cal TML}$ to show independent proofs of some (known) interesting properties of the logic ${\cal TML}$. In what follows  $\Gamma$, $\Delta$ are sets of formulas and $\alpha$, $\beta$, $\psi$ are formulas. \\ 
In the first place, we shall present a new independent proof of Proposition \ref{PropCompTML}. To do this, we need the following technical result. 

\begin{prop}\label{PropSecM4}  If $\vdash_{{\bf SC}_{\cal TML}} \Gamma \Rightarrow \Delta$ then,  for every
$\mathbb{A} \in {\bf TMA}$ and for every $h \in Hom(\mathfrak{Fm}, \mathbb{A})$, $\bigwedge_{\gamma \in \Gamma}  h(\gamma) \leq \bigvee_{\delta\in\Delta} h(\delta)$.
\end{prop}
\begin{dem} Suppose that   $\vdash_{{\bf SC}_{\cal TML}} \Gamma \Rightarrow \Delta$ and let $\cal P$ be a cut--free proof of the sequent $\Gamma \Rightarrow \Delta$ in ${\bf SC}_{\cal TML}$. Let $\mathbb{A} \in {\bf TMA}$ and let $h \in Hom(\mathfrak{Fm}, \mathbb{A})$. We use induction on the number $n$  of inferences in $\cal P$. If $n=0$ the proposition is obviously valid. (I.H.) Suppose that the proposition holds for $n<k$, $k>0$. Let $n=k$ and let $(r)$ be the last inference in $\cal P$. Ir $(r)$ is the right/left weakening rule, the proposition holds since $\mathbb{A}$ is, in particular, a lattice. If $(r)$ is one of the rules  {\rm ($\vee\Rightarrow$), ($\Rightarrow\vee$), ($\neg\vee\Rightarrow$), ($\Rightarrow\neg\vee$),  ($\wedge\Rightarrow$), ($\Rightarrow\wedge$), ($\neg\wedge\Rightarrow$), ($\Rightarrow\neg\wedge$),  ($\neg\neg\Rightarrow$), ($\Rightarrow\neg\neg$)}, the proposition holds since  $\mathbb{A}$ is, in particular, a De Morgan algebra. Finally, if $(r)$ is one of the rules {\rm , ($\square\Rightarrow$)$_i$, ($\Rightarrow\square$), ($\neg\square\Rightarrow$), ($\Rightarrow\neg\square$)$_i$ \, $i=1,2$} then the proposition holds since  $\mathbb{A}$  is a tetravalent modal algebra. For instance, suppose that $(r)$ is  $(\Rightarrow\square)$ and the last inference of $\cal P$ is $\displaystyle\frac{\Gamma\Rightarrow\Delta,\alpha \hspace{0.5cm}\Gamma, \neg\alpha\Rightarrow\Delta }{\Gamma \Rightarrow \Delta, \square\alpha}$. By (I.H.), we have (1) $\bigwedge_{\gamma \in \Gamma}h(\gamma) \leq \bigvee_{\delta\in\Delta}h(\delta) \vee h(\alpha)$ and  (2) $\bigwedge_{\gamma \in \Gamma}h(\gamma) \wedge h(\neg\alpha) \leq \bigvee_{\delta\in\Delta}h(\delta)$. Then, from (1), (2) and Proposition \ref{propsquare} we have  $\bigwedge_{\gamma \in \Gamma}h(\gamma)  \leq \bigvee_{\delta\in\Delta}h(\delta) \vee h(\square \alpha)$. 
\end{dem}

\

\begin{prop} The following conditions are equivalent.
\begin{itemize}
\item[(i)]  $\Gamma \models_{\cal TML} \psi$ ,
\item[(ii)] $\Gamma \models_{{\cal M}_4} \psi$.
\end{itemize}
\end{prop}
\begin{dem} {\em (i) imples (ii)}: immediate.\\
{\em (ii) implies (i)}: It is consequence of  Theorem \ref{TeoComSC} (i) and Proposition \ref{PropSecM4}.
\end{dem}

\

\noindent Next, we shall prove that the rule ($\neg$) of Font and Rius’ system is addmissible in ${\bf SC}_{\cal TML}$. Let $X$ a set of formulas, we shall denote by $\neg X$ the set  
$\neg X=\{\neg \gamma : \gamma \in X\}$. 

\

\begin{theo}\label{TheoContraPositive}  If $\vdash_{{\bf SC}_{\cal TML}} \Gamma \Rightarrow \Delta$, then $\vdash_{{\bf SC}_{\cal TML}} \neg\Delta\Rightarrow\neg\Gamma$.
\end{theo}
\begin{dem} Suppose that $\vdash_{{\bf SC}_{\cal TML}} \Gamma \Rightarrow \Delta$ and let $\cal P$ a cut--free proof of the sequent $\Gamma \Rightarrow \Delta$. We use induction on the number $n$ of inferences in ${\cal P}$. If $n=0$, then $\Gamma \Rightarrow \Delta$ is $\alpha\Rightarrow\alpha$, for some $\alpha$, and $\neg\Delta\Rightarrow\neg\Gamma$ is $\neg\alpha\Rightarrow\neg\alpha$ which is provable in ${\bf SC}_{\cal TML}$. (I.H.) Suppose that the lemma holds for $n<k$, with $k>0$. Let $n=k$ and let $(r)$ be the last inference in $\cal P$. If $(r)$ is left weakening, then the last inference of $\cal P$ is $\displaystyle\frac{\Gamma\Rightarrow\Delta}{\Gamma, \alpha \Rightarrow \Delta}$. By (I.H.), $\neg \Delta \Rightarrow \neg \Gamma$ is provable in ${\bf SC}_{\cal TML}$ and using right weakening we have $\vdash_{{\bf SC}_{\cal TML}} \neg\Delta\Rightarrow\neg\Gamma, \neg \alpha$. If $(r)$ is an instance of the right weakening the treatment is analogous.\\
Suppose now that $(r)$ is (an instance of) a logic rule. If $(r)$ is $(\Rightarrow\vee)$ and the last inference of $\cal P$ is  $\displaystyle\frac{\Gamma\Rightarrow\Delta, \alpha, \beta}{\Gamma \Rightarrow \Delta, \alpha\vee\beta}$. By (I.H.),  $\neg \alpha, \neg\beta, \neg\Delta \Rightarrow \neg \Gamma$ is provable,  and using $(\neg\vee\Rightarrow)$ we have that  $ \neg (\alpha\vee\beta), \neg\Delta \Rightarrow \neg \Gamma$ is provable. The cases where $(r)$ is one of the rules $(\vee\Rightarrow)$, $(\Rightarrow\neg\vee)$, $(\neg\vee\Rightarrow)$, $(\Rightarrow\wedge)$   $(\wedge\Rightarrow)$, $(\Rightarrow\neg\wedge)$, $(\neg\wedge\Rightarrow)$ are left to the reader.\\
If $(r)$ is $(\Rightarrow\neg\neg)$ and the last inference of $\cal P$ is  $\displaystyle\frac{\Gamma\Rightarrow\Delta, \alpha}{\Gamma \Rightarrow \Delta, \neg\neg\alpha}$. By (I.H.), $\neg\alpha, \neg\Delta\Rightarrow \neg\Gamma$ is provable in ${\bf SC}_{\cal TML}$ and using $(\neg\neg\Rightarrow)$ we have that  $\neg\neg\neg\alpha, \neg\Delta\Rightarrow \neg\Gamma$ is provable. If $(r)$ is  $(\neg\neg\Rightarrow)$ the proof is analogous.\\
If $(r)$ is $(\square\Rightarrow)_1$ and the last inference of $\cal P$ is $\displaystyle\frac{\Gamma, \alpha\Rightarrow\Delta}{\Gamma, \square\alpha \Rightarrow \Delta}$. By (I.H.), we have that  $\neg\Delta\Rightarrow \neg\Gamma, \neg\alpha$ is provable in ${\bf SC}_{\cal TML}$. Then, using $(\Rightarrow\neg\square)_2$ we have that  $\neg\Delta\Rightarrow \neg\Gamma, \neg\square\alpha$ is provable.\\
If $(r)$ is $(\square\Rightarrow)_2$ and the last inference of $\cal P$ is   $\displaystyle\frac{\Gamma\Rightarrow\Delta, \neg\alpha}{\Gamma, \square\alpha \Rightarrow \Delta}$.  By (I.H.), we have that  $\neg\Delta, \neg\neg\alpha\Rightarrow \neg\Gamma$ is provable in ${\bf SC}_{\cal TML}$ and using left weakening we have (1) $\vdash_{{\bf SC}_{\cal TML}} \alpha, \neg\Delta, \neg\neg\alpha\Rightarrow \neg\Gamma$. On the other hand, one can easily check that  $\vdash_{{\bf SC}_{\cal TML}} \alpha\Rightarrow \neg\neg\alpha$ and by means of (right/left) weakening(s) we have (2) $\vdash_{{\bf SC}_{\cal TML}} \alpha, \neg\Delta\Rightarrow \neg\neg\alpha, \neg\Gamma$. From (1), (2) and the cut rule, we have  $\vdash_{{\bf SC}_{\cal TML}} \alpha, \neg\Delta\Rightarrow \neg\Gamma$ (the cut rule is admissible in ${\bf SC}_{\cal TML}$). Then, using $(\Rightarrow\neg\square)_1$ we have  $\vdash_{{\bf SC}_{\cal TML}} \neg\Delta\Rightarrow \neg\Gamma, \neg\square\alpha$.\\
If  $(r)$ is $(\Rightarrow\square)$ and the last inference of $\cal P$ is $\displaystyle\frac{\Gamma\Rightarrow\Delta,\alpha \hspace{0.5cm}\Gamma, \neg\alpha\Rightarrow\Delta }{\Gamma \Rightarrow \Delta, \square\alpha}$. By (I.H.) we have that (3) $\vdash_{{\bf SC}_{\cal TML}} \neg\alpha, \neg\Delta\Rightarrow \neg\Gamma$ and (4) $\vdash_{{\bf SC}_{\cal TML}}  \neg\Delta\Rightarrow \neg\neg\alpha, \neg\Gamma$. From (4) and a similar reasoning to the above, we have that (5) $\vdash_{{\bf SC}_{\cal TML}}  \neg\Delta\Rightarrow \alpha, \neg\Gamma$. From (3), (5) and $(\neg\square\Rightarrow)$ we get $\vdash_{{\bf SC}_{\cal TML}}  \neg\Delta, \neg\square\alpha\Rightarrow \neg\Gamma$.\\
The cases where $(r)$ is one of the rules $(\Rightarrow\neg\square)_1$, $(\Rightarrow\neg\square)_2$ and $(\neg\square\Rightarrow)$ are treated similarly. 
\end{dem}

\

\begin{cor} $(\neg)$ is admissible in ${\bf SC}_{\cal TML}$.
\end{cor}

\noindent Finally,

\begin{theo}   $\vdash_{\cal TML}  \square\psi$ \, iff \,  $\vdash_{\cal TML}  \psi$.
\end{theo}
\begin{dem}($\Longrightarrow$) Suppose that $\vdash_{\cal TML}  \square\psi$. By Theorem \ref{TeoComSC}, Proposition \ref{PropCompTML} we know that the sequent $\Rightarrow \square\psi$ has a cut-free proof $\cal P$  in ${\bf SC}_{\cal TML}$. Let $(r)$ be the last inference of $\cal P$. By inspecting the rules of  ${\bf SC}_{\cal TML}$ we may assert that $(r)$ has to be an instance of the rule ($\Rightarrow\square$). So, $(r)$ is $\displaystyle \frac{\Rightarrow \psi \hspace{0.5cm} \neg\psi\Rightarrow}{\Rightarrow\square\psi}$ and clearly the sequent $\Rightarrow \psi$ is provable in ${\bf SC}_{\cal TML}$. Therefore $\vdash_{\cal TML}  \psi$.\\[2mm]
($\Longleftarrow$) Suppose that $\vdash_{\cal TML}  \psi$. By Theorem \ref{TeoComSC} (i), we have: (1) \, $\Rightarrow  \psi$ is provable in ${\bf SC}_{\cal TML}$. From (1) and  Theorem \ref{TheoContraPositive}, we have that: (2) \, $\neg\psi\Rightarrow$ is also provable in ${\bf SC}_{\cal TML}$. From (1), (2) and the rule ($\Rightarrow\square)$,  we may assert that $\Rightarrow \square\psi$ is provable in  ${\bf SC}_{\cal TML}$. Therefore,  $\vdash_{\cal TML}  \square\psi$.
\end{dem}

\section{Natural deduction for ${\cal TML}$}

In this section, we shall present a natural deduction system for ${\cal TML}$. We take our inspiration from the construction  made before. In particular, it threw some light on how the connective $\square$ behaves. We think that this system shows an interesting example of a rule (different from the usual ones), namely the introduction rule of the connective $\square$,  that needs to produce a discharge of hypothesis; and this is related to the intrinsic meaning of the connective. \\
The proof system ${\bf ND}_{\cal TML}$ will be defined following the notational conventions given in \cite{TS}.

\begin{defi} Deductions in ${\bf ND}_{\cal TML}$ are inductively defined as follows:\\
Basis: The proof tree with a single occurrence of an assumption $\phi$ with a marker is a deduction with conclusion $\phi$ from open assumption $\phi$ .\\[2mm]
Inductive step: Let ${\cal D}$, ${\cal D}_1$ ,${\cal D}_2$,${\cal D}_3$ be deductions. Then, they can be extended by one of the following rules below. The classes {\rm[$\neg\phi$]$^u$}, {\rm[$\neg\psi$]$^v$}, {\rm[$\phi$]$^u$} , {\rm[$\psi$]$^v$} below contain open assumptions of the deductions of the premises of the final inference, but are closed in the whole deduction. 

\begin{prooftree}
\AxiomC{}
\RightLabel{\rm MA (modal axioma)}
\UnaryInfC{$\phi \vee \neg \square\phi$}
\end{prooftree}

\begin{prooftree}
\AxiomC{${\cal D}_1$}
\noLine
\UnaryInfC{$\phi$}

\AxiomC{${\cal D}_2$}
\noLine
\UnaryInfC{$\psi$}
\RightLabel{\rm $\wedge$I}
\BinaryInfC{$\phi \wedge \psi$}

\AxiomC{}

\AxiomC{${\cal D}$}
\noLine
\UnaryInfC{$\phi\wedge\psi$}
\RightLabel{\rm $\wedge$E$_1$}
\UnaryInfC{$\phi$}

\AxiomC{}

\AxiomC{${\cal D}$}
\noLine
\UnaryInfC{$\phi\wedge\psi$}
\RightLabel{\rm $\wedge$E$_2$}
\UnaryInfC{$\psi$}

\noLine
\QuinaryInfC{}
\end{prooftree}

\begin{prooftree}
\AxiomC{${\cal D}$}
\noLine
\UnaryInfC{$\neg\phi$}
\RightLabel{\rm $\neg\wedge$I$_1$}
\UnaryInfC{$\neg(\phi\wedge\psi)$}

\AxiomC{}

\AxiomC{${\cal D}$}
\noLine
\UnaryInfC{$\neg\psi$}
\RightLabel{\rm $\neg\wedge$I$_2$}
\UnaryInfC{$\neg(\phi\wedge\psi)$}

\AxiomC{}

\AxiomC{}
\noLine
\UnaryInfC{${\cal D}_1$}
\noLine
\UnaryInfC{$\neg(\phi\wedge\psi)$}

\AxiomC{\rm[$\neg\phi$]$^u$}
\noLine
\UnaryInfC{${\cal D}_2$}
\noLine
\UnaryInfC{$\chi$}

\AxiomC{\rm [$\neg\psi$]$^v$}
\noLine
\UnaryInfC{${\cal D}_3$}
\noLine
\UnaryInfC{$\chi$}
\RightLabel{\rm $\neg\wedge$E,$u$,$v$}
\TrinaryInfC{$\chi$}
\noLine
\QuinaryInfC{}
\end{prooftree}

\begin{prooftree}

\AxiomC{${\cal D}$}
\noLine
\UnaryInfC{$\phi$}
\RightLabel{\rm $\vee$I$_1$}
\UnaryInfC{$\phi\vee\psi$}

\AxiomC{}

\AxiomC{${\cal D}$}
\noLine
\UnaryInfC{$\psi$}
\RightLabel{\rm $\vee$I$_2$}
\UnaryInfC{$\phi\vee\psi$}

\AxiomC{}

\AxiomC{}
\noLine
\UnaryInfC{${\cal D}_1$}
\noLine
\UnaryInfC{$\phi\vee\psi$}

\AxiomC{\rm[$\phi$]$^u$}
\noLine
\UnaryInfC{${\cal D}_2$}
\noLine
\UnaryInfC{$\chi$}

\AxiomC{\rm [$\psi$]$^v$}
\noLine
\UnaryInfC{${\cal D}_3$}
\noLine
\UnaryInfC{$\chi$}
\RightLabel{\rm $\vee$E,$u$,$v$}
\TrinaryInfC{$\chi$}
\noLine
\QuinaryInfC{}
\end{prooftree}

\begin{prooftree}

\AxiomC{${\cal D}_1$}
\noLine
\UnaryInfC{$\neg\phi$}

\AxiomC{${\cal D}_2$}
\noLine
\UnaryInfC{$\neg\psi$}
\RightLabel{\rm $\neg\vee$I}
\BinaryInfC{$\neg(\phi \vee \psi)$}

\AxiomC{}

\AxiomC{${\cal D}$}
\noLine
\UnaryInfC{$\neg(\phi\vee\psi)$}
\RightLabel{\rm $\neg\vee$E$_1$}
\UnaryInfC{$\neg\phi$}

\AxiomC{}

\AxiomC{${\cal D}$}
\noLine
\UnaryInfC{$\neg(\phi\vee\psi)$}
\RightLabel{\rm $\neg\vee$E$_2$}
\UnaryInfC{$\neg\psi$}

\noLine
\QuinaryInfC{}
\end{prooftree}

\begin{prooftree}

\AxiomC{${\cal D}$}
\noLine
\UnaryInfC{$\phi$}
\RightLabel{\rm $\neg\neg$I}
\UnaryInfC{$\neg\neg\phi$}

\AxiomC{}

\AxiomC{}

\AxiomC{}

\AxiomC{${\cal D}$}
\noLine
\UnaryInfC{$\neg\neg\phi$}
\RightLabel{\rm $\neg\neg$E}
\UnaryInfC{$\phi$}

\noLine
\QuinaryInfC{}
\end{prooftree}

\begin{prooftree}

\AxiomC{${\cal D}_1$}
\noLine
\UnaryInfC{$\psi\vee\phi$}

\AxiomC{\rm [$\neg\phi$]$^u$}
\noLine
\UnaryInfC{${\cal D}_2$}
\noLine
\UnaryInfC{$\psi$}
\RightLabel{\rm $\square$I$^*$,$u$}
\BinaryInfC{$\psi\vee\square\phi$}

\AxiomC{}

\AxiomC{}

\AxiomC{}

\AxiomC{${\cal D}$}
\noLine
\UnaryInfC{$\square\phi$}
\RightLabel{\rm $\square$E}
\UnaryInfC{$\phi$}

\noLine
\QuinaryInfC{}
\end{prooftree}

\begin{prooftree}

\AxiomC{${\cal D}$}
\noLine
\UnaryInfC{$\neg\phi$}
\RightLabel{\rm $\neg\square$I}
\UnaryInfC{$\neg\square\phi$}

\AxiomC{}

\AxiomC{}

\AxiomC{}

\AxiomC{${\cal D}_1$}
\noLine
\UnaryInfC{$\neg\square\phi$}

\AxiomC{${\cal D}_2$}
\noLine

\UnaryInfC{$\phi$}

\RightLabel{\rm $\neg\square$E}
\BinaryInfC{$\neg\phi$}

\noLine
\QuinaryInfC{}
\end{prooftree}

\begin{prooftree}

\AxiomC{${\cal D}$}
\noLine
\UnaryInfC{$\neg\phi\wedge\square\phi $}
\RightLabel{$\bot$I}
\UnaryInfC{$\bot$}

\AxiomC{}

\AxiomC{}

\AxiomC{}

\AxiomC{${\cal D}$}
\noLine
\UnaryInfC{$\bot $}
\RightLabel{$\bot$E}
\UnaryInfC{$\alpha$}
\noLine
\QuinaryInfC{}
\end{prooftree}

\end{defi}

\

\begin{rem} If we take $\psi$ as $\bot$ in  $\square$I$^*$ we get  \begin{prooftree}
\AxiomC{${\cal D}_1$}
\noLine
\UnaryInfC{$\phi$}

\AxiomC{\rm [$\neg\phi$]$^u$}
\noLine
\UnaryInfC{${\cal D}_2$}
\noLine
\UnaryInfC{$\bot$}
\RightLabel{\rm $\square$I,$u$}
\BinaryInfC{$\square\phi$}
\end{prooftree}
\noindent Formally,  $\square$I is derivable in ${\bf ND}_{\cal TML}$.  The intuition behind this rule is the following:``if we have a deduction for $\alpha$ and $\neg \alpha$ is not provable, then we have a deduction for $\square\alpha$''. \\
\end{rem}

As usual, by application of the rule $\neg\wedge$E a new proof-tree is formed from ${\cal D}$, ${\cal D}_1$, and ${\cal D}_2$ by adding at the bottom the conclusion $\chi$ while closing the sets {\rm [$\neg\phi$]$^u$} and {\rm [$\neg\psi$]$^u$} of open  assumptions  marked by $u$ and $v$, respectively. Idem for the rules $\wedge$E and $\square$I. Note that we have introduced the symbol $\bot$, it behaves here as an arbitrary unprovable propositional constant. \\[2mm]
Let $\Gamma \cup \{\alpha\} \subseteq Fm$. We say that the conclusion $\alpha$ is derivable from a set $\Gamma$ of premises, noted $\Gamma\vdash \alpha$, if and only if there is a deduction in ${\bf ND}_{\cal TML}$ of $\alpha$ from $\Gamma$.

\

\begin{theo}(Soundness and Completeness) Let $\Gamma, \Delta\subseteq Fm$, $\Gamma$ finite. The following conditions are equivalent:
\begin{itemize} 
\item[(i)] the sequent \ $\Gamma\Rightarrow\Delta$ \ is derivable in ${\bf SC}_{\cal TML}$,
\item[(ii)] there is a deduction of the disjunction of the sentences in $\Delta$ from $\Gamma$ in ${\bf ND}_{\cal TML}$.
\end{itemize}
\end{theo}
\begin{dem} (i) implies (ii): Suppose that the sequent $\Gamma\Rightarrow\Delta$ is derivable in ${\bf SC}_{\cal TML}$, that is, there is a formal proof $\cal P$ of $\Gamma\Rightarrow\Delta$ in ${\bf SC}_{\cal TML}$ which does not use the cut rule. We shall show that there is a deduction  of the disjunction of the formulas in $\Delta$ (denoted by $\bigvee\Delta$) from $\Gamma$ in ${\bf ND}_{\cal TML}$,  using induction on the number $n$ of rule applications in $\cal P$, $n\geq 0$. \\
If $n=0$, then $\Gamma\Rightarrow\Delta$ is $\alpha\Rightarrow\alpha$ and it is clear that $\alpha\vdash\alpha$. Now, (I.H.) suppose that ``(i) implies (ii)'' holds for $n<k$, with $k>0$.\\
Let $n=k$, that is $\cal P$ is a derivation in ${\bf SC}_{\cal TML}$ with last rule (r) of the form
\begin{prooftree}
\AxiomC{$\Gamma_1\Rightarrow\Delta_1$}
\noLine
\UnaryInfC{$\ddots$}

\AxiomC{$\dots$}
\noLine
\UnaryInfC{$\vdots$}

\AxiomC{$\Gamma_t\Rightarrow\Delta_t$}
\noLine
\UnaryInfC{\reflectbox{$\ddots$}}

\noLine
\TrinaryInfC{}
\noLine
\UnaryInfC{}
\LeftLabel{\small ($r$)}
\UnaryInfC{$\Gamma\Rightarrow\Delta$}
\end{prooftree}
\noindent If (r) is left weakening, then the last rule of $\cal P$ has the form $\displaystyle {\rm (r)}\frac{\Gamma'\Rightarrow\Delta}{\Gamma', \beta\Rightarrow\Delta}$. By (I.H.), there exists a deduction $\cal D$ of $\Delta$ from $\Gamma'$, then 
\begin{prooftree}
\AxiomC{$\cal D$}
\noLine
\UnaryInfC{$\bigvee\Delta$}
\AxiomC{$\beta$}
\RightLabel{\small $\wedge$I}
\BinaryInfC{$\bigvee\Delta\wedge\beta$}
\RightLabel{\small $\wedge$E$_1$}
\UnaryInfC{$\bigvee\Delta$}
\end{prooftree}
\noindent is a deduction of $\bigvee\Delta$ from $\Gamma'\cup\{\beta\}$. If (r) is right weakening, then (r) has the form $\displaystyle {\rm (r)}\frac{\Gamma\Rightarrow\Delta'}{\Gamma\Rightarrow\Delta', \beta}$, then by (I.H.) there is a deduction $\cal D$ of $\bigvee\Delta'$ from $\Gamma$.
\begin{prooftree}
\AxiomC{$\cal D$}
\noLine
\UnaryInfC{$\bigvee\Delta'$}
\RightLabel{\small $\vee$I$_1$}
\UnaryInfC{$\bigvee\Delta'\vee\beta$}
\end{prooftree}
Now, suppose that (r) is a logical rule, we shall prove it just for ($\Rightarrow\vee$), ($\vee\Rightarrow$), ($\Rightarrow\neg\vee$), ($\neg\vee\Rightarrow$). If (r) is ($\vee\Rightarrow$), then we may assume that the last inference of $\cal P$ has the form  $\displaystyle \mbox{($\Rightarrow\vee$)}\frac{\Gamma\Rightarrow \Delta', \alpha, \beta}{\Gamma\Rightarrow\Delta', \alpha \vee \beta}$. Then, by (I.H.) we have a deduction $\cal D$ of $\bigvee \Delta' \vee \alpha \vee \beta$ from $\Gamma$ and the proof is complete.\\
If (r) is ($\vee\Rightarrow$) and last inference of $\cal P$ has the from $\displaystyle \mbox{($\vee\Rightarrow$)}\frac{\Gamma, \gamma_1\Rightarrow \Delta \hspace{.5cm}\Gamma,\gamma_2\Rightarrow \Delta}{\Gamma, \gamma_1 \vee \gamma_2 \Rightarrow \Delta}$, then by (I.H.) there are deductions ${\cal D}_i$ , $i=1,2$, of $\alpha$ from $\Gamma\cup \{\gamma_i\}$. Then, the following
\begin{prooftree}
\AxiomC{$\gamma_1 \vee \gamma_2$}

\AxiomC{[$\gamma_1$]$^{u_1}$}
\noLine
\UnaryInfC{${\cal D}_1$}
\noLine
\UnaryInfC{$\bigvee\Delta$}

\AxiomC{[$\gamma_2$]$^{u_2}$}
\noLine
\UnaryInfC{${\cal D}_2$}
\noLine
\UnaryInfC{$\bigvee \Delta$}
\RightLabel{\small $\vee$E,$u_1$,$u_2$}
\TrinaryInfC{$\bigvee\Delta$}
\end{prooftree}
is a deduction of $\bigvee\Delta$ from $\Gamma\cup\{\gamma_1\vee\gamma_2\}$. Note that in this last deduction we have made every assumption $\gamma_i$ in ${\cal D}_i$ an open assumption with label $u_i$.\\
If (r) is ($\neg\vee\Rightarrow$) then we may assume that the last instance of ${\cal P}$ has the form $\displaystyle \mbox{($\neg\vee\Rightarrow$)}\frac{\Gamma, \neg\gamma_1\Rightarrow \Delta }{\Gamma, \neg(\gamma_1 \vee \gamma_2) \Rightarrow \Delta}$. By (I.H.), there is a deduction ${\cal D}$ of $\alpha$ from $\Gamma\cup\{\gamma_1\}$ and the following
\begin{prooftree}
\AxiomC{$\neg(\gamma_1 \vee \gamma_2)$}
\RightLabel{\small $\neg\vee$E$_1$}
\UnaryInfC{$\neg\gamma_1$}
\noLine
\UnaryInfC{${\cal D}$}
\noLine
\UnaryInfC{$\bigvee\Delta$}
\end{prooftree}
is a deduction of $\alpha$ from $\Gamma\cup\{\neg(\gamma_1 \vee \gamma_2)\}$. If (r) is  ($\Rightarrow\neg\vee$) we proceed analogously.\\
For (r) being any of the rules ($\square\Rightarrow$)$_i$, ($\Rightarrow\square$), ($\neg\square\Rightarrow$), ($\Rightarrow\neg\square$)$_i$ \, $i=1,2$, we present the next table showing the deduction corresponding to the premise(s) of (r) and the deduction corresponding to the consequence of (r).    

\

{\small \hspace{-1cm}
\begin{tabular}{|l|c|c|}  \hline
&  & \\
 Rule (r)    &  Upper sequent(s)'s    &  Lower sequent's  \\
                & deduction(s)    & deduction \\
& &  \\ \hline \hline
& &  \\ 
{\small \AxiomC{$\Gamma, \gamma \Rightarrow \Delta$}
\LeftLabel{\tiny ($\square\Rightarrow$)$_1$}
\UnaryInfC{$\Gamma, \square\gamma \Rightarrow \Delta$}
\DisplayProof}
 &
\AxiomC{$\gamma$}
\noLine
\UnaryInfC{$\cal D$}
\noLine
\UnaryInfC{$\bigvee\Delta$}
\DisplayProof
 & {\tiny
\AxiomC{$\square\gamma$}
\RightLabel{\tiny $\square$E}
\UnaryInfC{$\gamma$}
\noLine
\UnaryInfC{$\cal D$}
\noLine
\UnaryInfC{$\bigvee \Delta$}
\DisplayProof} \\

& &  \\ \hline
& & \\
{\small\AxiomC{$\Gamma \Rightarrow \Delta, \neg\gamma$}
\LeftLabel{\tiny ($\square\Rightarrow$)$_2$}
\UnaryInfC{$\Gamma, \square\gamma \Rightarrow \Delta$}
\DisplayProof}
 &
\AxiomC{$\cal D$}
\noLine
\UnaryInfC{$\bigvee\Delta\vee\neg\gamma$}
\DisplayProof
& {\tiny
\AxiomC{$\cal D$}
\noLine
\UnaryInfC{$\bigvee\Delta\vee\neg\gamma$}

\AxiomC{$[\bigvee\Delta]_u$}
\noLine
\UnaryInfC{$\bigvee\Delta$}

\AxiomC{$[\neg\gamma]_u$}
\AxiomC{$\square\gamma$}
\RightLabel{\tiny $\wedge$I}
\BinaryInfC{$\neg\gamma\wedge\square\gamma$}
\RightLabel{\tiny $\bot$I}
\UnaryInfC{$\bot$}
\RightLabel{\tiny $\bot$E}
\UnaryInfC{$\bigvee \Delta$}
\RightLabel{$\vee$E, $u$, $v$}
\TrinaryInfC{$\bigvee \Delta$}
\DisplayProof}\\

& &  \\ \hline
& & \\
{\small \AxiomC{$\Gamma \Rightarrow \Delta, \gamma$}
\AxiomC{$\Gamma,\neg \gamma \Rightarrow \Delta $}
\LeftLabel{\tiny ($\Rightarrow\square$)}
\BinaryInfC{$\Gamma \Rightarrow \Delta, \square\gamma$}
\DisplayProof}
&

\AxiomC{${\cal D}_1$}
\noLine
\UnaryInfC{$\bigvee\Delta\vee\gamma$}
\DisplayProof
 \hskip 1.5em
\AxiomC{$\neg\gamma$}
\noLine
\UnaryInfC{${\cal D}_2$}
\noLine
\UnaryInfC{$\bigvee \Delta$}
\DisplayProof

& 
{\tiny
\AxiomC{${\cal D}_1$}
\noLine
\UnaryInfC{$\bigvee\Delta\vee\gamma$}

\AxiomC{$\neg\gamma^u$}
\noLine
\UnaryInfC{${\cal D}_2$}
\noLine
\UnaryInfC{$\bigvee\Delta$}
\RightLabel{\tiny $\square$I$^*$,$u$}
\BinaryInfC{$\bigvee\Delta\vee\square\gamma$}
\DisplayProof }\\
& &  \\ \hline
& & \\
{\small \AxiomC{$\Gamma \Rightarrow \Delta,\gamma$}
\AxiomC{$\Gamma,\neg \gamma \Rightarrow \Delta$}
\LeftLabel{\tiny ($\neg\square\Rightarrow$)}
\BinaryInfC{$\Gamma,\neg\square\gamma \Rightarrow\Delta$}
\DisplayProof}
&
\AxiomC{${\cal D}_1$}
\noLine
\UnaryInfC{$\bigvee\Delta\vee\gamma$}
\DisplayProof

\AxiomC{$\neg\gamma$}
\noLine
\UnaryInfC{${\cal D}_2$}
\noLine
\UnaryInfC{$\bigvee\Delta$}
\DisplayProof
& {\tiny
\AxiomC{${\cal D}_1$}
\noLine
\UnaryInfC{$\bigvee\Delta\vee\gamma$}

\AxiomC{$[\neg\gamma]^u$}
\noLine
\UnaryInfC{${\cal D}_2$}
\noLine
\UnaryInfC{$\bigvee\Delta$}
\RightLabel{\tiny $\square$I,$u$}
\BinaryInfC{$\bigvee\Delta\vee\square\gamma$}

\AxiomC{$\neg\square \gamma$}
\RightLabel{\tiny $\vee$I$_2$}
\UnaryInfC{$\bigvee\Delta\vee\neg\square\gamma$}
\RightLabel{\tiny $\wedge$I}
\BinaryInfC{$(\bigvee\Delta\vee\square\gamma)\wedge(\bigvee\Delta\vee\neg\square\gamma)$}
\RightLabel{\footnotemark}
\doubleLine
\UnaryInfC{$\bigvee\Delta\vee(\square\gamma\wedge\neg\square\gamma)$}
\RightLabel{\footnotemark}
\doubleLine
\UnaryInfC{$\bigvee\Delta\vee\bot$}
\UnaryInfC{$\bigvee\Delta$}

\DisplayProof} \\

& &  \\ \hline
& & \\
{\small \AxiomC{$\Gamma,\gamma \Rightarrow \Delta$}
\LeftLabel{\tiny ($\Rightarrow\neg\square$)$_1$}
\UnaryInfC{$\Gamma\Rightarrow\Delta, \neg\square\gamma$}
\DisplayProof}
&
\AxiomC{$\gamma$}
\noLine
\UnaryInfC{$\cal D$}
\noLine
\UnaryInfC{$\bigvee\Delta$}
\DisplayProof
&
{\tiny
\AxiomC{}
\RightLabel{\tiny (MA)}
\UnaryInfC{$\gamma\vee\neg\square\gamma$}

\AxiomC{$\gamma^v$}
\noLine
\UnaryInfC{$\cal D$}
\noLine
\UnaryInfC{$\bigvee\Delta$}
\UnaryInfC{$\bigvee\Delta\vee\neg\square\gamma$}

\AxiomC{$\neg\square \gamma^u$}
\UnaryInfC{$\bigvee\Delta\vee\neg\square\gamma$}

\RightLabel{\tiny $\vee$E,$u$,$v$}
\TrinaryInfC{$\bigvee\Delta\vee\neg\square\gamma$}
\DisplayProof} \\

& &  \\ \hline
& & \\
{\small \AxiomC{$\Gamma \Rightarrow\Delta, \neg\gamma$}
\LeftLabel{\tiny ($\Rightarrow\neg\square$)$_2$}
\UnaryInfC{$\Gamma\Rightarrow\Delta, \neg\square\gamma$}
\DisplayProof}
&
\AxiomC{$\cal D$}
\noLine
\UnaryInfC{$\bigvee\Delta\vee\neg\gamma$}
\DisplayProof
& 
{\tiny
\AxiomC{$\cal D$}
\noLine
\UnaryInfC{$\bigvee\Delta\vee\neg\gamma$}

\AxiomC{$\bigvee\Delta^u$}
\UnaryInfC{$\bigvee\Delta\vee\neg\square\gamma$}

\AxiomC{$\neg\gamma^v$}
\RightLabel{\tiny $\neg\square$I}
\UnaryInfC{$\neg\square\gamma$}
\UnaryInfC{$\bigvee\Delta\vee\neg\square\gamma$}

\RightLabel{\tiny $\vee$E,$u$,$v$}
\TrinaryInfC{$\bigvee\Delta\vee\neg\square\gamma$}

\DisplayProof}  \\ 

 & & \\ \hline

\end{tabular} 
}

\

\noindent (ii) implies (i):  Let $\cal D$ be a deduction of the disjunction of the sentences in $\Delta$ from $\Gamma$ in ${\bf ND}_{\cal TML}$.  As before, we use induction on the number $n$ of rule instances in the deduction $\cal D$. If $r=0$ the proof is trivial. (I.H.) Suppose that ``(ii) implies (i)'' holds for $n<k$, $k>0$; and let $(r)$ the last rule instance in $\cal D$. If $(r)$ is one of the introduction/elimination rule of $\wedge$I, $\wedge$E, $\neg\wedge$I, $\neg\wedge$E, $\vee$I, $\vee$E , $\neg\vee$I, $\neg\vee$E, $\neg\neg$I and $\neg\neg$E; the proof is immediate since these rules are just translations of the corresponding rules of ${\bf SC}_{\cal TML}$. Suppose that $(r)$ is $\square$I$^*$, then $\cal D$ is

\begin{prooftree}
\AxiomC{${\cal D}_1$}
\noLine
\UnaryInfC{$\psi\vee\phi$}

\AxiomC{\rm [$\neg\phi$]$^u$}
\noLine
\UnaryInfC{${\cal D}_2$}
\noLine
\UnaryInfC{$\psi$}
\RightLabel{\rm $\square$I$^*$,$u$}
\BinaryInfC{$\psi\vee\square\phi$}
\end{prooftree}

\noindent Then, by (I.H), we have that the sequents $\Gamma_1\Rightarrow \psi\vee\phi$ and $\Gamma_2, \neg\phi \Rightarrow \psi$ are provable in   ${\bf SC}_{\cal TML}$,  where $\Gamma_1\cup\Gamma_2 =\Gamma$. By using weakening(s) and the cut rule we obtain $\Gamma\Rightarrow \psi, \phi$ and $\Gamma \neg\phi \Rightarrow \psi$ are provable. Then, using ($\square\Rightarrow$), we have that $\vdash_{{\bf SC}_{\cal TML}} \Gamma\Rightarrow\psi, \square\phi$.  If $(r)$ is $\square$E, then $\cal D$ is

\begin{prooftree}
\AxiomC{${\cal D}$}
\noLine
\UnaryInfC{$\square\phi$}
\RightLabel{\rm $\square$E}
\UnaryInfC{$\phi$}
\end{prooftree}

\noindent By (I.H.), we have $\vdash_{{\bf SC}_{\cal TML}} \Gamma\Rightarrow \square\phi$. From the fact that $\vdash_{{\bf SC}_{\cal TML}}  \square\phi\Rightarrow \phi$ and the cut rule the proof is completed.  If $(r)$ is $\neg\square$I, then $\cal D$ is
\begin{prooftree}
\AxiomC{${\cal D}$}
\noLine
\UnaryInfC{$\neg\phi$}
\RightLabel{\rm $\neg\square$I}
\UnaryInfC{$\neg\square\phi$}
\end{prooftree}

\noindent By (I.H.), we have $\vdash_{{\bf SC}_{\cal TML}} \Gamma\Rightarrow \neg\phi$. By Theorem \ref{TheoContraPositive}, $\vdash_{{\bf SC}_{\cal TML}} \neg\neg\phi\Rightarrow \neg\Gamma$ and from $\vdash_{{\bf SC}_{\cal TML}} \phi\Rightarrow  \neg\neg\phi$ and the cut rule, we have $\vdash_{{\bf SC}_{\cal TML}} \phi\Rightarrow \neg\Gamma$. Using ($\square\Rightarrow)$ we obtain  $\vdash_{{\bf SC}_{\cal TML}} \square\phi\Rightarrow \neg\Gamma$ and by Theorem \ref{TheoContraPositive}  $\vdash_{{\bf SC}_{\cal TML}} \neg\neg\Gamma\Rightarrow \neg\square\phi$. Finally, from  $\vdash_{{\bf SC}_{\cal TML}} \Gamma\Rightarrow \neg\neg\Gamma$ and cut(s) (and weakening(s) if necessary) we obtain $\vdash_{{\bf SC}_{\cal TML}}\Gamma\Rightarrow \neg\square\phi$. If $(r)$ is $\neg\square$E, then $\cal D$ is
\begin{prooftree}
\AxiomC{${\cal D}_1$}
\noLine
\UnaryInfC{$\neg\square\phi$}

\AxiomC{${\cal D}_2$}
\noLine

\UnaryInfC{$\phi$}

\RightLabel{\rm $\neg\square$E}
\BinaryInfC{$\neg\phi$}
\end{prooftree}

\noindent  By (I.H) and using weakening(s) we have that the sequents $\Gamma\Rightarrow \neg\square\phi$ and $\Gamma \Rightarrow \phi$ are provable in   ${\bf SC}_{\cal TML}$. Using ($\Rightarrow\wedge$), we obtain $\vdash_{{\bf SC}_{\cal TML}} \Gamma\Rightarrow  \phi\wedge\neg\square\phi$ and since $\phi\wedge\neg\square\phi \Leftrightarrow \phi\wedge\neg\phi$ and the cut rule we obtain  $\vdash_{{\bf SC}_{\cal TML}} \Gamma\Rightarrow  \phi\wedge\neg\phi$. Finally, taking into account that  $\vdash_{{\bf SC}_{\cal TML}} \phi\wedge\neg\phi\Rightarrow \neg\phi$ we have $\vdash_{{\bf SC}_{\cal TML}} \Gamma\Rightarrow\neg\phi$.\\
The cases in which $(r)$ is $\bot$I or $\bot$E are immediate (see Remark \ref{RemBot}).
\end{dem}

\footnotetext[1]{$(\gamma \vee\alpha)\wedge(\gamma\vee\beta) \dashv\vdash \gamma\vee (\alpha \wedge  \beta)$}
 \footnotetext[2]{ $\alpha \vee (\square\gamma \wedge \neg\square\gamma) \dashv\vdash \alpha\vee \bot$}
\

Since our natural deduction system is strongly inspired by the cut-free sequent calculus ${\bf SC}_{\cal TML}$, one can likely expect normalization to hold for ${\bf SC}_{\cal TML}$. 

\section{Conclusions}

In the present paper we focused on the  proof-theoretic aspects of the tetravalent modal logic ${\cal TML}$. In the first place, we showed that the strongly adequate Gentzen calculus given by Font and Rius for ${\cal TML}$ does not enjoy the cut--elimination property. Then, by applying a method due to Avron, Ben-Naim and Konikowska, we developed a sequent calculus for $\cal TML$ with the cut--elimination property. This allowed us to provide new independent proof of some known interesting properties of ${\cal TML}$. Finally, strongly inspired by this cut--free sequent calculus,  we presented a natural deduction system, sound and complete with respect to the ${\cal TML}$.

Despite the fact that ${\cal TML}$ was originally defined as the logic that preserves degrees of truth w.r.t. tetravalent modal algebras, we could use  Avron, Ben-Naim and Konikowska's method; and this is because  ${\cal TML}$ is also a matrix logic.
 An interesting task to be done is to extend this method to logics to logics  that preserves degrees of truth w.r.t. some ordered structure but which do not have a matrix semantics.

\section{Acknowledgments}  I would like to thank the anonymous referees for their extremely careful reading, helpful suggestions and constructive comments on this paper. 
\

\end{document}